\documentclass[11pt]{article}
\usepackage[english]{babel}
\usepackage{lineno}
\usepackage{graphicx,multicol}
\usepackage{epic,eepic,epsfig}
\usepackage{caption}

\usepackage{latexsym,bm}
\usepackage{mathrsfs}
\usepackage{amsmath,amsthm}
\usepackage{graphicx}
\usepackage{amssymb}
\usepackage{CJK}
\usepackage{color}
\usepackage{cite}
\usepackage{comment}
\usepackage{epsfig}
\usepackage{epstopdf}
\usepackage{algorithm}
\usepackage{algorithmicx}
\usepackage{algpseudocode}
\usepackage{slashed}

\floatname{algorithm}{Algorithm}

\usepackage{tikz}

\pgfdeclarelayer{edgelayer}
\pgfdeclarelayer{nodelayer}
\pgfsetlayers{edgelayer,nodelayer,main}

\tikzstyle{none}=[inner sep=0pt]
\definecolor{hexcolor0xf81e1c}{rgb}{0.973,0.118,0.110}
\definecolor{hexcolor0x3c00ff}{rgb}{0.235,0.000,1.000}

\tikzstyle{whitevertex}=[circle,fill=white,draw=black, scale = 0.5]
\tikzstyle{redvertex}=[circle,fill=hexcolor0xf81e1c,draw=black, scale = 0.5]
\tikzstyle{bluevertex}=[circle,fill=hexcolor0x3c00ff,draw=black, scale = 0.5]
\tikzstyle{greenvertex}=[circle,fill=green,draw=black, scale=0.5]
\tikzstyle{purplevertex}=[circle,fill=magenta,draw=black, scale=0.5]
\tikzstyle{grayvertex}=[circle,fill=white,draw=gray, scale=0.5]
\tikzstyle{blackvertex}=[circle,fill=black,draw=black, scale=0.5]

\tikzstyle{textbox}=[rectangle,fill=none,draw=none]
\tikzstyle{box}=[rectangle,fill=none,draw=black]

\tikzstyle{arc}=[black, ->]
\tikzstyle{grayarc}=[gray, ->]
\tikzstyle{bluearc}=[blue, ->]
\tikzstyle{grayedge}=[draw=gray]
\tikzstyle{blueedge}=[draw=blue]
\tikzstyle{rededge}=[draw=red]
\tikzstyle{edge}=[draw=black]

\tikzstyle{vertex}=[circle, ,fill=white,draw=black, scale=0.5]

\tikzstyle{10circle}=[circle, scale=10.0,draw=black]
\tikzstyle{10oval}=[ellipse, scale=10.0,draw=black]

\theoremstyle{plain}
\newtheorem{thm}{\hspace{5mm}Theorem}[section]
\newtheorem{cor}[thm]{\hspace{5mm}Corollary}
\newtheorem{lem}[thm]{\hspace{5mm}Lemma}
\newtheorem{prop}[thm]{\hspace{5mm}Proposition}
\newtheorem{prob}[thm]{\hspace{5mm}Problem}

\newcommand{\pf}{{\bf Proof: }}

\theoremstyle{definition}

\theoremstyle{plain}

\theoremstyle{plain}

\theoremstyle{plain}

\theoremstyle{plain}

\makeatletter
\newenvironment{breakablealgorithm}
  {
   \begin{center}
     \refstepcounter{algorithm}
     \hrule height.8pt depth0pt \kern2pt
     \renewcommand{\caption}[2][\relax]{
       {\raggedright\textbf{\ALG@name~\thealgorithm} ##2\par}%
       \ifx\relax##1\relax 
         \addcontentsline{loa}{algorithm}{\protect\numberline{\thealgorithm}##2}%
       \else 
         \addcontentsline{loa}{algorithm}{\protect\numberline{\thealgorithm}##1}%
       \fi
       \kern2pt\hrule\kern2pt
     }
  }{
     \kern2pt\hrule\relax
   \end{center}
  }
\makeatother
\DeclareGraphicsExtensions{.eps,.eps.gz}
\normalsize \rm
\begin{document}
\begin{CJK}{GBK}{song}
\newcommand{\song}{\CJKfamily{song}}    
\newcommand{\fs}{\CJKfamily{fs}}        
\newcommand{\kai}{\CJKfamily{kai}}      
\newcommand{\hei}{\CJKfamily{hei}}      
\newcommand{\li}{\CJKfamily{li}}        

\title{Comparability digraphs: An analogue of comparability graphs
\footnote{This work is supported by NSFC (Grant No. 12071194) and NSERC.}}
\author{Xiao-Lu Gao\thanks{School of Mathematics and Statistics, 
                           Lanzhou University,
                           Lanzhou, Gansu 730000, China; gaoxl14@lzu.edu.cn}, 
       Jing Huang\thanks{Department of Mathematics and Statistics, 
                         University of Victoria, 
                         Victoria, B.C. V8W 2Y2, Canada; huangj@uvic.ca},
       Shou-Jun Xu\thanks{School of Mathematics and Statistics, 
                          Lanzhou University,
                          Lanzhou, Gansu 730000, China; shjxu@lzu.edu.cn} }

\date{}

\maketitle
\thispagestyle{empty}

\begin{abstract}
Comparability graphs are a popular class of graphs. We introduce as the digraph
analogue of comparability graphs the class of comparability digraphs. We show that 
many concepts such as implication classes and the knotting graph for a comparability
graph can be naturally extended to a comparability digraph. We give 
a characterization of comparability digraphs in terms of their knotting graphs. 
Semicomplete comparability digraphs are a prototype of comparability digraphs. 
One instrumental technique for analyzing the structure of comparability graphs is 
the Triangle Lemma for graphs. We generalize the Triangle Lemma to semicomplete 
digraphs. Using the Triangle Lemma for semicomplete digraphs we prove that if 
an implication class of a semicomplete digraph contains no circuit of length 2 then 
it contains no circuit at all. We also use it to device an $\mathcal{O}(n^3)$ 
time recognition algorithm for semicomplete comparability digraphs where $n$ 
is the number of vertices of the input digraph. The correctness of the algorithm 
implies a characterization for semicomplete comparability digraphs, akin to that 
for comparability graphs.

\end{abstract}

\section{Introduction}

Many graph classes such as interval graphs, circular-arc graphs, chordal graphs, 
split graphs, threshold graphs, and cographs have found their corresponding digraph 
counterparts: 
interval digraphs \cite{sdrw}, circular-arc digraphs \cite{Circular-arc-di}, 
chordal digraphs \cite{chordal-di}, split digraphs \cite{lamar}, oriented threshold 
graphs \cite{boe}, and directed cographs \cite{bdr}. 
These digraph counterparts not only generalize the corresponding graphs but also 
carry much rich structural properties for the study of many graph theory problems.
In this paper, we introduce the digraph counterpart of comparability graphs.

Comparability graphs are the graphs which have transitive orientations. These graphs
can be equivalently defined in terms of vertex orderings. 
A graph $G = (V,E)$ is a {\em comparability graph} if it has a vertex ordering 
$\prec$ such that for all $x \prec y \prec z$, $xy \in E$ and $yz \in E$ imply
$xz \in E$. Such an ordering $\prec$ is called a {\em comparability ordering} of $G$.
When $G$ has a comparability ordering $\prec$, orienting each edge $xy$ from 
$x$ to $y$ whenever $x \prec y$ gives a transitive orientation of $G$. 

Let $D = (V,A)$ be a digraph. We call $D$ a {\em comparability digraph} if 
it has a vertex ordering $\prec$ such that for all $x \prec y \prec z$, 
$xy \in A$ and $yz \in A$ imply $xz \in A$ and $zy \in A$ and $yx \in A$ imply 
$zx \in A$. We call $\prec$ a {\em comparability ordering} of $D$. 
An arc $xy$ of $D$ is {\em symmetric} if $yx$ is also an arc of $D$. If every arc
of $D$ is symmetric then $D$ is a {\em symmetric} digraph. Graphs are essentially 
symmetric digraphs. With this view every comparability graph is a comparability 
digraph. Thus comparability digraphs are a generalization of comparability graphs.

The naturalness of comparability digraphs as a generalization of comparability 
graphs can also be seen from their adjacency matrices. Let $D$ be a digraph and 
$M(D)$ be an adjacency matrix of $D$. Note that the entries of $M(D)$ on the main 
diagonal are all 0. We say that the two-by-two identity matrix 
$I_2 = \left[\begin{matrix}1&0\\
                           0&1
                \end{matrix}\right]$  
is a {\em principal submatrix} of $M(D)$ if either of the entries 0 of $I_2$ lies on 
the main diagonal. It turns out that a comparability ordering of $D$ corresponds to 
a simultaneous permutation of rows and columns of $M(D)$ so that the resulting form 
of the matrix does not contain $I_2$ as a principal submatrix. In particular,
a graph $G$ is a comparability graph if and only if its adjacency matrix $M(G)$ 
admits a simultaneous permutation of rows and columns so that $I_2$ is not 
a principal submatrix of the resulting matrix. In general, $0,1$-matrices which 
admit permutations of rows and columns so that the resulting matrices do not contain
$I_2$ as a submatrix have been characterized in \cite{hhlim}. Bipartite graphs whose 
biadjacency matrices having this property are called {\em comparability bigraphs} 
\cite{hhlim}. Matrices which avoid other small forbidden submatrices via permutations
of rows and columns have been studied in 
\cite{afarber,dein,minorder,hks,krw,lubiw1}. 

Every transitive digraph is a comparability digraph but not all acyclic digraphs 
are comparability digraphs (see example in Figure~\ref{noC2}). A tournament is 
a comparability digraph if and only if it is acyclic. Every semicomplete 
comparability digraph contains an acyclic tournament as a spanning subdigraph.
If a digraph is a comparability digraph then the subdigraph induced by the symmetric
arcs is a comparability digraph. Every strong digraph which contains no odd directed
cycle has a bipartite underlying graph and hence is a comparability digraph.

Comparability graphs have elaborate structural properties. Many of them have been 
exploited by Gallai in his seminal paper \cite{gallai}.
Gallai \cite{gallai} proved that a graph is a comparability graph if and only if 
its {\em knotting graph} (see definition in Section~\ref{ik}) is bipartite. 
This implies that comparability graphs can be recognized in polynomial time. 
Faster recognition algorithms for comparability graphs can be found in 
\cite{golumbic77,recg-lexico,lss}. A linear time algorithm which finds a transitive
orientation of a comparability graph is deviced in \cite{1999Modular}.
A forbidden subgraph characterization of comparability graphs has been given
in \cite{gallai}.

A key concept in the study of comparability graphs is the notion of implication 
classes. The implication classes of a graph $G$ form a partition of the set of
ordered pairs of adjacent vertices. The knotting graph of $G$ is such a graph whose 
components consist of implication classes of $G$. If the knotting graph of $G$ is 
not bipartite then $G$ is not a comparability graph. On the other hand, if 
the knotting graph of $G$ is bipartite, then each component of the knotting graph
corresponds to a maximal subgraph of $G$ that has an essentially unique transitive 
orientation. A transitive orientation of $G$ is made of suitable transitive
orientation of the components of the knotting graph. Thus, comparability graphs
can be characterized in terms of implication classes and knotting graphs.
We extend the concepts of implication classes and knotting graphs of graphs to 
digraphs. We show comparability digraphs can also be characterized in terms of 
knotting graphs (see Theorem \ref{knotting-char}). Unfortunately it appears 
bewildered that this characterization leads to a polynomial time recognition 
algorithm for comparability digraphs. It seems challenging to find such 
an algorithm.  

Semicomplete digraphs are a special class of digraphs which admit beautiful 
properties as well as polynomial time solutions to many problems 
hard for general digraphs. Given a graph $G$ and a vertex ordering $\prec$ of
$G$, let $D$ be the semicomplete digraph obtained from $G$ by replacing each edge 
$xy$ of $G$ with the symmetric arcs $xy,yx$ and each non-edge $xy$ of $G$ with 
non-symmetric arc $xy$ whenever $x \prec y$. It is easy to see that if $\prec$
is a comparability ordering of $G$ then $\prec$ is a comparability ordering of $D$. 
On the other hand, if $\prec$ is a comparability ordering of a semicomplete digraph 
$D$, then $\prec$ is a comparability ordering of the graph obtained $D$ by replacing
each symmetric arc $xy$ of $D$ with an edge $xy$. Thus, as far as the property of 
being comparability is concerned, semicomplete digraphs overpower graphs.

The Triangle Lemma for graphs is useful for devicing an efficient recognition 
algorithm for comparability graphs (cf. \cite{golumbic}). We show that 
the Triangle Lemma for semicomplete digraphs (see Lemma~\ref{lem:Triangle}) 
is equally useful for the investigation of implication classes of semicomplete 
digraphs and is intrumental for the device of an efficient recognition algorithm 
for semicomplete comparability digraphs.  

The rest of the paper is organized as follows. In Section~\ref{ik}, we extend 
the concepts implication classes and knotting graphs useful for graphs to digraphs. 
We characterize comparability digraphs in terms of knotting graphs. A necessary 
condition for a digraph to be a comparability digraph is that none of its 
implication classes contains a circuit. We show by example this condition is not
sufficient. We extend, in Section~\ref{semicomplete}, the Triangle Lemma for graphs 
to semicomplete digraphs. Using the Triangle Lemma we prove that if an implication 
class of a semicomplete digraph contains no circuit of length 2 then it contains no 
circuit at all. In Section \ref{recog}, we device an $\mathcal{O}(n^3)$ time 
recognition algorithm for semicomplete comparability digraphs where $n$ is the number
of vertices of the input digraph. The correctness of the algorithm allows us to
derive a list of characterizations of semicomplete comparability digraphs similar to 
those for comparability graphs. Finally, we conclude in Section~\ref{conclusion} 
with some remarks and open problems.

\section{Implication classes and the knotting graph}
\label{ik}

Implication classes and knotting graphs are important for analyzing the 
structure of comparability graphs, cf. \cite{gallai,golumbic77,golumbic}. 
We extend these concepts to digraphs.  

For a digraph $D=(V,A)$, let $Z_D$ denote the set of ordered pairs of adjacent 
vertices in $D$, that is,
\[Z_D = \{(x,y):\ xy \in A\ \mbox{or}\ yx \in A\}. \]
Note that each arc $xy$ of $D$ gives rise to two ordered pairs $(x,y)$ and $(y,x)$ 
in $Z_D$. For simplicity we call each ordered pair in $Z_D$ a {\em pair}. 
Given pairs $(x,y), (x',y')$ in $Z_D$, we say that $(x,y)$ {\em directly forces} 
$(x',y')$ and denote it by $(x,y) \Gamma (x',y')$ if one of the following holds:

\begin{itemize}
\item $x = x'$ and $y = y'$;
\item $x = x'$, $y \neq y'$, and either $yx, x'y' \in A$ and $yy' \notin A$ or 
                          $y'x', xy \in A$ and $y'y \notin A$;
\item $y = y'$, $x \neq x'$, and either $xy, y'x' \in A$ and $xx' \notin A$ or
                          $x'y', yx \in A$ and $x'x \notin A$.
\end{itemize} 

Clearly, $\Gamma$ is a reflexive and symmetric relation on $Z_D$. The transitive 
closure $\Gamma^*$ of $\Gamma$ is thus an equivalence relation on $Z_D$. Each 
equivalence class of $Z_D$ under $\Gamma^*$ is called an {\em implication class} of 
$D$. We use $\mathcal{I}(D)$ to denote the set of implication classes of $D$.
Two pairs $(x,y)$ and $(x',y')$ are in the same implication class (i.e.,
$(x,y) \Gamma^* (x',y')$) if and only if there exists a sequence of pairs 
$(x_0,y_0), (x_1,y_1), \dots, (x_k,y_k)$ such that
\[(x,y) = (x_0,y_0) \Gamma (x_1,y_1) \Gamma \cdots \Gamma (x_k,y_k) = (x',y').\]
Such a sequence is called a {\em $\Gamma$-chain} connecting $(x,y)$ and $(x',y')$.
All these definitions apply to graphs $G$ which are viewed as symmetric digraphs, 
in which case $Z_G$ and $\mathcal{I}(G)$ are used for the sets of pairs and 
implication classes of $G$ respectively.  

A $\Gamma$-chain $(x_0,y_0), (x_1,y_1), \dots, (x_k,y_k)$ is called {\em canonical} 
if $x_i = x_{i+1}$ and $y_i = y_{i-1}$ for all odd $i$.

\begin{prop} \label{canon}
If $(x,y) \Gamma^* (x',y')$, then there is a canonical $\Gamma$-chain connecting 
$(x,y)$ and $(x',y')$. 
\end{prop}
\pf Let $(x_0,y_0), (x_1,y_1), \dots, (x_k,y_k)$ be a $\Gamma$-chain connecting
$(x,y)$ and $(x',y')$. Obtain a new chain from this chain by adding $(x_{i+1},y_i)$ 
between $(x_i,y_i)$ and $(x_{i+1},y_{i+1})$ for each $0\leq i < k$. That is, 
$(x_{i+1},y_i)$ immediately succeeds $(x_i,y_i)$ and immediately precedes 
$(x_{i+1},y_{i+1})$ in the new chain. Since $(x_i,y_i) \Gamma (x_{i+1},y_{i+1})$, 
$x_i = x_{i+1}$ or $y_i = y_{i+1}$ for each $i$. Thus $(x_{i+1},y_i)$ is equal to 
$(x_i,y_i)$ or $(x_{i+1},y_{i+1})$ for each $i$ and hence
$(x_i,y_i) \Gamma (x_{i+1},y_i) \Gamma (x_{i+1},y_{i+1})$. Moreover, each added pair
has its first vertex equal to the first vertex of the pair immediately succeeding it
and its second vertex equal to the second vertex of the pair immediately preceeding 
it. Therefore the new chain is by definition a canonical $\Gamma$-chain connecting 
$(x,y)$ and $(x',y')$.
\qed

In case when $D$ is a comparability digraph, pairs in the same implication class of 
$D$ are consistent with respect to each comparability ordering of $D$. 
  
\begin{prop} \label{force}
Let $D$ be a comparability digraph and $\prec$ a comparability ordering of $D$. 
Suppose $(x,y) \Gamma^* (x',y')$. Then $x \prec y$ if and only if $x' \prec y'$.
\end{prop} 
\pf It suffices to show that $x \prec y$ if and only if $x' \prec y'$ for pairs
$(x,y), (x',y')$ with $(x,y) \Gamma (x',y')$. This is clearly true when $x = x'$ and
$y = y'$. When $x = x'$ and $y \neq y'$, either $yx, x'y' \in A$ and 
$yy' \notin A$ or $y'x', xy \in A$ and $y'y \notin A$. Since $\prec$ is 
a comparability ordering of $D$, $x \prec y$ if and only if $x' \prec y'$. 
A similar proof applies when $y = y'$ and $x \neq x'$.  
\qed

Let $I \subseteq Z_D$ be a set of pairs. The {\em inverse} of $I$, denoted by 
$I^{-1}$, is the following set
\[\{(y,x):\ (x,y)\in I\}. \] 
Clearly, $I \in \mathcal{I}(D)$ implies $I^{-1} \in \mathcal{I}(D)$, i.e.,
the inverse of an implication class is an implication class.
A vertex ordering $\prec$ of $D$ may be viewed as the set of pairs 
$(x,y) \in Z_D$ with $x \prec y$. In this way Proposition~\ref{force} simply
asserts that each comparability ordering of $D$ is a union of implication classes 
that contains one of $I, I^{-1}$ for each $I \in \mathcal{I}(D)$.  
A {\em circuit} of {\em length $k$} is a set of pairs 
$(x_1,x_2), (x_2,x_3), \dots, (x_{k-1},x_k), (x_k,x_1)$.
When the pairs of a circuit are all contained in $I$, we say that 
$I$ {\em contains} the circuit.

\begin{prop} \label{ob:NoCirInvIntEmp}
For any $I \in \mathcal{I}(D)$, $I = I^{-1}$ if and only if $I$ contains a circuit 
of length 2.
\end{prop}
\pf If $I = I^{-1}$, then for any $(x,y) \in I$, $(y,x) \in I$ so $(x,y), (y,x)$ is 
a circuit of length 2 contained in $I$. Conversely, suppose $(x,y),(y,x)$ is 
a circuit of length 2 contained in $I$. Then, for each $(x',y') \in I$,
$(x,y)\Gamma^* (x',y')\Gamma^* (y,x)$, which implies $(y',x') \in I$. 
Hence $I = I^{-1}$.
\qed

\begin{prop}\label{pro:NoCirComImp}
Suppose that $D$ is a comparability digraph. Then no implication class of $D$ 
contains a circuit. In particular, $I\neq I^{-1}$ for all $I \in \mathcal{I}(D)$.
\end{prop}
\pf Let $\prec$ be a comparability ordering of $D$ and let $I \in \mathcal{I}(D)$. 
If $(x_1,x_2), (x_2,x_3), \dots, (x_k,x_1)$ is a circuit contained in $I$, then 
\[(x_1,x_2) \Gamma^* (x_2,x_3) \Gamma^* \cdots \Gamma^* (x_k,x_1).\]
By Proposition \ref{force}, either
\[x_1 \prec x_2 \prec x_3 \prec \cdots \prec x_k \prec x_1\]
or
\[x_1 \succ x_2 \succ x_3 \succ \cdots \succ x_k \succ x_1.\]
In either case we have $x_1 \prec x_1$, which is a contradiction.
Hence $I$ contains no circuit and by Proposition \ref{ob:NoCirInvIntEmp},
$I\neq I^{-1}$.
\qed

\begin{figure}[ht]
\begin{center}
\begin{tikzpicture}[>=latex]
                \node [style=blackvertex] (1) at (-2,1) {};
                \node [style=blackvertex] (2) at (-1,1) {};
                \node [label={above:$w$}] [style=blackvertex] (3) at (0,1) {};
                \node [style=blackvertex] (4) at (1,1) {};
                \node [style=blackvertex] (5) at (2,1) {};
                \node [label={left:$x$}] [style=blackvertex] (6) at (-1,0) {};
                \node [label={right:$z$}] [style=blackvertex] (7) at (1,0) {};
                \node [style=blackvertex] (8) at (-2,-1) {};
                \node [style=blackvertex] (9) at (-1,-1) {};
                \node [label={below:$y$}] [style=blackvertex] (10) at (0,-1) {};
                \node [style=blackvertex] (11) at (1,-1) {};
                \node [style=blackvertex] (12) at (2,-1) {};

                \draw [-> ] (1) -> (2);
                \draw [-> ] (2) to (3);
                \draw [-> ] (3) to (4);
                \draw [-> ] (4) to (5);
                \draw [-> ] (2) to (6);
                \draw [-> ] (3) to (6);
                \draw [-> ] (3) to (7);
                \draw [-> ] (7) to (4);
                \draw [->] (6) to (9);
                \draw [->] (10) to (6);
                \draw [->] (10) to (7);
                \draw [->] (11) to (7);
                \draw [->] (12) to (11);
                \draw [->] (11) to (10);
                \draw [->] (10) to (9);
                \draw [->] (9) to (8);
\end{tikzpicture}
\end{center}
\vspace{-2mm}
\caption{A non-comparability digraph $D$ for which $I\neq I^{-1}$ for all
$I \in \mathcal{I}(D)$.}
\label{noC2}
\end{figure}
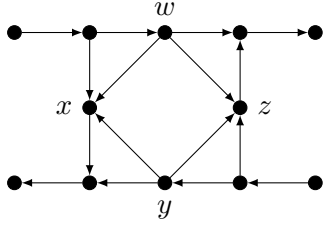

The digraph in Figure~\ref{noC2} has an implication class which contains circuit 
$(w,x), (x,y), (y,z), (z,w)$ so it is not a comparability digraph according to
Proposition~\ref{pro:NoCirComImp}.

\begin{thm} \cite{golumbic} \label{conv}
A graph $G$ is a comparability graph if and only if $I\neq I^{-1}$ for all
$I \in \mathcal{I}(G)$.
\qed
\end{thm}

The non-comparability digraph $D$ in Figure~\ref{noC2} satisfies the property that 
$I \neq I^{-1}$ for all $I \in \mathcal{I}(D)$. So Theorem \ref{conv} is not true 
for digraphs. 

We now turn to the notion of knotting graphs. Let $D$ be a digraph and let $x$ be 
a vertex of $D$. Define a relation on the set of pairs $(x,y)$ of $Z_D$ at $x$ as 
follows: Two pairs $(x,y)$ and $(x,z)$ at $x$ are {\em knotted}, denoted by 
$(x,y) \tilde{\Gamma} (x,z)$, if there exists a sequence of pairs
$(x,y_0), (x,y_1), \dots, (x,y_k)$ at $x$ such that
\[(x,y) = (x,y_0) \Gamma (x,y_1) \Gamma \cdots \Gamma (x,y_k) = (x,z).\]
Note that $(x,y) \tilde{\Gamma} (x,z)$ implies $(x,y) \Gamma^* (x,z)$ but the 
converse is not true in general. Nevertheless, it is easy to see that 
$\tilde{\Gamma}$ is an equivalence relation on the set of pairs at $x$. 
This equivalence relation yields naturally a partition of the neighbourhood $N(x)$ 
of $x$ into sets $x^1, x^2, \dots, x^{\ell_x}$ in such a way that for all 
$y,z \in N(x)$, $(x,y) \tilde{\Gamma}(x,z)$ if and only if $y, z \in x^{\alpha}$ 
for some $1 \leq \alpha \leq \ell_x$. This partition of $N(x)$ is well-defined for 
each vertex $x$ of $D$. Moreover, if $x, y$ are adjacent vertices of $D$, then 
there exist unique $\alpha, \beta$ with $1 \leq \alpha \leq \ell_x$ and 
$1 \leq \beta \leq \ell_y$ such that $x \in y^{\beta}$ and $y \in x^{\alpha}$.

\begin{figure}[ht]
\begin{center}
\begin{tikzpicture}[>=latex]
%

                \node [style=blackvertex] (13) at (3.8,1.2) {};
                \node [style=blackvertex] (14) at (5,1.2) {};
                \node [label={above:$w^1$}] [style=blackvertex] (15) at (6.2,1.2) {};
                \node [style=blackvertex] (16) at (7.4,1.2) {};
                \node [style=blackvertex] (17) at (8.6,1.2) {};
                \node [label={below:$w^2$}] [style=blackvertex] (18) at (6.2,.9) {};

                \node [label={left:$x^1$}] [style=blackvertex] (19) at (5,0) {};
                \node [label={right:$x^2$}] [style=blackvertex] (20) at (5.3,0) {};
                \node [label={left:$z^1$}] [style=blackvertex] (21) at (7.1,0) {};
                \node [label={right:$z^2$}] [style=blackvertex] (22) at (7.4,0) {};

                \node [label={above:$y^1$}] [style=blackvertex] (23) at (6.2,-.9) {};
                \node [style=blackvertex] (24) at (3.8,-1.2) {};
                \node [style=blackvertex] (25) at (5,-1.2) {};
                \node [label={below:$y^2$}] [style=blackvertex] (26) at (6.2,-1.2) {};
                \node [style=blackvertex] (27) at (7.4,-1.2) {};
                \node [style=blackvertex] (28) at (8.6,-1.2) {};

               \draw [-] (13) -- (14) -- (15) -- (16) -- (17);
               \draw [-] (14) -- (19) -- (25);
               \draw [-] (16) -- (22) -- (27);
               \draw [-] (24) -- (25) -- (26) -- (27) -- (28);
               \draw [-] (18) -- (19);
               \draw [-] (20) -- (26);
               \draw [-] (22) -- (23);
               \draw [-] (15) -- (21);

\end{tikzpicture}
\end{center}
\vspace{-2mm}
\caption{The knotting graph of the digraph in Figure~\ref{noC2}.}
\label{knotting}
\end{figure}
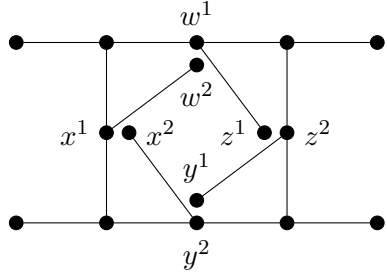

The {\em knotting graph} $\tilde{K}_D$ of $D$ is the graph having vertices
$x^1, x^2, \dots, x^{\ell_x}$ for all $x \in V(D)$ and edges $x^{\alpha}y^{\beta}$ 
for all $x, y \in V(D)$ with $x \in y^{\beta}$ and $y \in x^{\alpha}$. 
We call $x^1, x^2, \dots, x^{\ell_x}$ {\em partial copies} of $x$ and call
$x$ the {\em original} of their partial copies. Thus each vertex $x$ of $D$ gives 
rise to $\ell_x$ vertices and each pair of adjacent vertices $x, y$ in $D$
gives rise to the unique edge $x^{\alpha}y^{\beta}$ in $\tilde{K}_D$. 
If $x^{\alpha}y^{\beta}, x^{\alpha}z^{\gamma}$ are two edges incident with 
$x^{\alpha}$ in $\tilde{K}_D$, then $(x,y) \tilde{\Gamma} (x,z)$. By identifying 
the partial copies of $x$ in $\tilde{K}_D$ for all $x \in V(D)$ we obtain 
the underlying graph $U(D)$ of $D$. The knotting graph of the digraph 
in Figure~\ref{noC2} is depicted in Figure~\ref{knotting}.

A {\em knotting walk} in $D$ is a sequence of vertices $x_0, x_1, \dots, x_p$ 
such that $(x_i,x_{i-1})$ and $(x_i,x_{i+1})$ are knotted for all $0 < i < p$. 
Each knotting walk $x_0, x_1, \dots, x_p$ corresponds naturally to the walk 
$x_0^{\alpha_0}, x_1^{\alpha_1}, \dots, x_p^{\alpha_p}$ in $\tilde{K}_D$ where 
each $x_i^{\alpha_i}$ is a partial copy of $x_i$, and $x_0 \in x_1^{\alpha_1}$,
$x_p \in x_{p-1}^{\alpha_{p-1}}$, and 
$x_i \in x_{i-1}^{\alpha_{i-1}} \cap x_{i+1}^{\alpha_{i+1}}$ for each 
$0 < i < p$. 
Conversely, if $x_0^{\alpha_0}, x_1^{\alpha_1}, \dots, x_p^{\alpha_p}$ is a walk 
in $\tilde{K}_D$ then $x_0, x_1, \dots, x_p$ is a knotting walk in $D$.
Knotting walks are another way to connecting pairs in the same implication class of
$D$. Indeed, if $x_0, x_1, \dots, x_p$ is a knotting walk in $D$, then
$(x_0,x_1) \Gamma^* (x_{p-1},x_p)$ when $p$ is odd, and 
$(x_0,x_1) \Gamma^* (x_p, x_{p-1})$ when $p$ is even. Any two pairs in the same 
implication class of $D$ are connected by a knotting walk.   

\begin{prop} \label{knotting-walk}
If $(x,y)$ and $(x',y')$ are in the same implication class of $D$, then there is 
a knotting walk $x_0, x_1, \dots, x_p$ in $D$ with $(x_0,x_1) = (x,y)$ and 
$(x_{p-1},x_p) = (x',y')$.
\end{prop}
\pf Since $(x,y)$ and $(x',y')$ are in the same implication class in $D$, 
$(x,y) \Gamma^* (x',y')$. By Proposition \ref{canon} there is a canonical
$\Gamma$-chain $(x_0,y_0), (x_1,y_1), \dots, (x_k,y_k)$ connecting $(x,y)$ and 
$(x',y')$. Since the chain is canonical, $x_i = x_{i+1}$ and $y_i = y_{i-1}$ for 
each odd $i$. Hence $x_0, y_0, x_1, y_1, \dots, x_k, y_k$ is a knotting walk in $D$ 
connecting $(x,y)$ and $(x',y')$. 
\qed

Each knotting walk in $D$ corresponds to a walk in $\tilde{K}_D$ which connects
two edges in the same component of $\tilde{K}_D$. So the edges in a component
of $\tilde{K}_D$ correspond to the pairs in $I \cup I^{-1}$ for some 
$I \in \mathcal{I}(D)$. It turns out $I \neq I^{-1}$ if and only if 
the corresponding component of $\tilde{K}_D$ is bipartite. This can be seen from 
the proof of Theorem~\ref{thm:necessity} below. 

\begin{thm}\label{thm:necessity}
The following statements are equivalent for a digraph $D$:
\begin{enumerate}
\item ${\tilde{K}}_{D}$ is a bipartite graph;
\item $D$ contains no odd closed knotting walk;
\item $I\neq I^{-1}$ for all $I\in\mathcal{I}(D)$.
\end{enumerate}
\end{thm}
\pf Each closed walk in $\tilde{K}_{D}$ corresponds to a closed knotting walk in $D$
of the same length. The knotting graph $\tilde{K}_{D}$ is bipartite if and only if
it contains no odd closed walk, that is, $D$ contains no odd closed knotting walk.
This shows statements $1$ and $2$ are equivalent.

Suppose that $I = I^{-1}$ for some $I\in\mathcal{I}(D)$. Let $(x,y)$ be a pair 
in $I$. Then $(y,x)$ is also in $I$. By Proposition \ref{knotting-walk},
there is a knotting walk $x_0, x_1, \dots, x_p$ in $D$ with $(x_0,x_1) = (x,y)$ and
$(x_{p-1},x_p) = (y,x)$. Thus $x_0 = x = x_p$ and $x_1 = y = x_{p-1}$, which
imply $p$ is even and the knotting walk is closed. Hence $D$ contains an odd closed 
knotting walk. On the other hand, suppose that $D$ contains an odd closed knotting 
walk. Let $x_0, x_1, ..., x_{2k}, x_0$ be such a walk.
Then $(x_i, x_{i-1}) \Gamma^* (x_i, x_{i+1})$ for each $1\leq i\leq 2k$ where
the subscripts are modulo $2k+1$. Hence
\[(x_0, x_1)\Gamma^* (x_2, x_1)\Gamma^* (x_2, x_3)\Gamma^* \cdots \Gamma^* 
(x_{2k}, x_{2k-1})\Gamma^* (x_{2k}, x_0) \Gamma^* (x_1,x_0).\]
Thus the implication class $I$ that contains $(x_0, x_1)$ also contains 
$(x_1,x_0)$, which means $I = I^{-1}$. Therefore statements $2$ and $3$ are also 
equivalent.
\qed

Combining Proposition \ref{pro:NoCirComImp} and Theorem \ref{thm:necessity} we have
the following:

\begin{cor}\label{cor:neccesity}
The statements in Theorem \ref{thm:necessity} hold for every comparability
digraph $D$.
\qed
\end{cor}

\begin{figure}[ht]
\begin{center}
\begin{tikzpicture}[>=latex]
                \node [label={left:$w$}] [style=blackvertex] (1) at (-5,0) {};
                \node [label={right:$v$}] [style=blackvertex] (2) at (7,0) {};
                \node [label={below:$u$}] [style=blackvertex] (3) at (0,.5) {};
                \node [style=blackvertex] (4) at (2,1) {};
                \node [style=blackvertex] (5) at (1,1.2) {};
                \node [style=blackvertex] (6) at (-1,2) {};
                \node [style=blackvertex] (7) at (1,2) {};
                \node [style=blackvertex] (8) at (0,3) {};
                \node [label={right:$z$}] [style=blackvertex] (9) at (1,3) {};
                \node [label={right:$y$}] [style=blackvertex] (10) at (.5,4) {};
                \node [label={left:$x$}] [style=blackvertex] (11) at (0,5) {};
                \node [style=blackvertex] (12) at (0,6) {};

                \draw [-> ] (1) -> (11);
                \draw [-> ] (2) to (1);
                \draw [-> ] (3) to (1);
                \draw [-> ] (6) to (1);
                \draw [-> ] (8) to (1);
                \draw [-> ] (10) to (1);
                \draw [-> ] (2) to (3);
                \draw [-> ] (2) to (4);
                \draw [-> ] (2) to (8);
                \draw [-> ] (2) to (11);
                \draw [-> ] (7) to (2);
                \draw [-> ] (10) to (2);
                \draw [-> ] (3) to (4);
                \draw [-> ] (6) to (3);
                \draw [-> ] (7) to (3);
                \draw [-> ] (8) to (3);
                \draw [-> ] (4) to (5);
                \draw [-> ] (2) to (4);
                \draw [-> ] (7) to (4);
                \draw [->] (7) to (5);
                \draw [->] (6) to (7);
                \draw [->] (7) to (8);
                \draw [->] (8) to (9);
                \draw [->] (10) to (9);
                \draw [->] (8) to (11);
                \draw [->] (10) to (8);
                \draw [->] (10) to (11);
                \draw [->] (11) to (12);
\end{tikzpicture}
\end{center}
\vspace{-2mm}
\caption{A non-comparability digraph $D$ for which no $I \in \mathcal{I}(D)$
        contains a circuit.}
\label{noCs}
\end{figure}
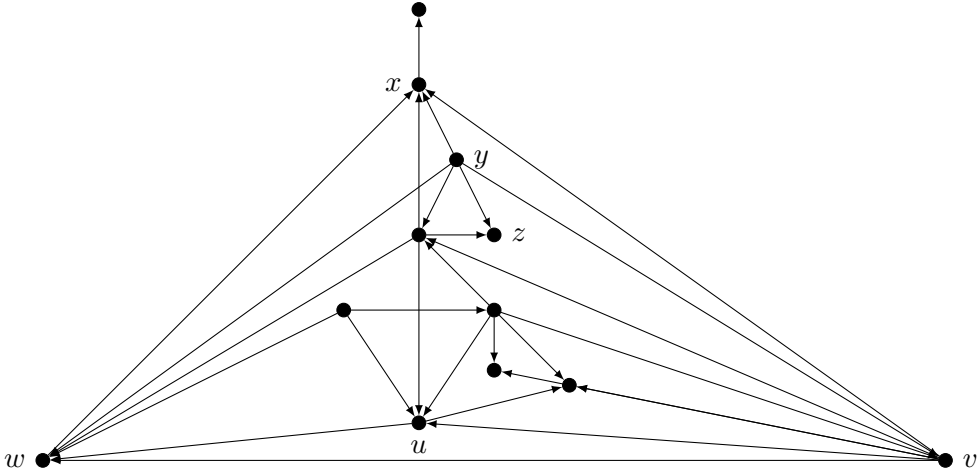

There are non-comparability digraphs for which none of their implication classes 
contains circuits so in particular they satisfy the statements in 
Theorem~\ref{thm:necessity}. The digraph in Figure~\ref{noCs} is such an example, 
in which $[x,y], [y,x], [u,v], [v,u]$ are the nontrivial implication classes and 
$ [y,w], [w,y], [y,z], [z,y]$ are the trivial implication classes. Although none of 
the implication classes contains a circuit, each of $[x,y]\cup [u,v]$ and 
$[x,y]\cup [v,u]$ contains one. So the digraph is not a comparability digraph. 

Surprisingly, the statements in Theorem~\ref{thm:necessity} are sufficient 
conditions for a graph to be a comparability graph.

\begin{thm} \cite{gallai, golumbic}\label{thm:undirected}
The following statements are equivalent for a graph $G$:
\begin{enumerate}
\item $G$ is a comparability graph;
\item the knotting graph $\tilde{K}_G$ is bipartite;
\item $G$ contains no odd closed knotting walk;
\item $I \neq I^{-1}$ for all $I \in \mathcal{I}(G)$;
\item $I$ contains no circuit for all $I \in \mathcal{I}(G)$.
\qed
\end{enumerate}
\end{thm}

Let $D$ be a digraph. Each edge $xy$ in $U(D)$ corresponds to an edge 
$x^{\alpha}y^{\beta}$ in $\tilde{K}_{D}$ where $x^{\alpha}$ and $y^{\beta}$ are 
partial copies of $x$ and $y$ respectively. This is a one-to-one correspondence 
between the edge set of $U(D)$ and the edge set of $\tilde{K}_{D}$.
An orientation of the edge $xy$ from $x$ to $y$ transforms to an orientation of 
the edge $x^{\alpha}y^{\beta}$ from $x^{\alpha}$ to $y^{\beta}$. Thus each 
orientation of $U(D)$ {\em transforms} into an orientation of $\tilde{K}_{D}$. 
This transformation is reversible. Indeed, given an orientation of $\tilde{K}_{D}$,
by identifying the partial copies of each vertex of $D$ we obtain an orientation
of $U(D)$.  

\begin{thm} \label{knotting-char}
A digraph $D$ is a comparability digraph if and only if $\tilde{K}_{D}$ is bipartite
and there exists a bipartition $(X,Y)$ of $\tilde{K}_{D}$ such that the orientation 
of $\tilde{K}_{D}$ from $X$ to $Y$ transforms into an acyclic orientation of $U(D)$.
\end{thm}
\pf Assume that $D$ is a comparability digraph. Let $\prec$ be a comparability 
ordering of $D$. Orient each edge $x^{\alpha}y^{\beta}$ of $\tilde{K}_{D}$ 
from $x^{\alpha}$ to $y^{\beta}$ if $x \prec y$ where $x$ and $y$ are the originals 
of $x^{\alpha}$ and $y^{\beta}$ respectively. This gives an orientation of 
$\tilde{K}_{D}$. Let $x^{\alpha}$ be a vertex and $y^{\beta}, z^{\gamma}$ be 
neighbours of $x^{\alpha}$ in $\tilde{K}_{D}$. Then $(x,y) \Gamma^* (x,z)$ where 
$x, y, z$ are the originals of $x^{\alpha}, y^{\beta}, z^{\gamma}$ respectively. 
By Proposition~\ref{force}, either $x \prec y$ and $x \prec z$ or 
$y \prec x$ and $z \prec x$. Hence in the orientation of $\tilde{K}_{D}$,
$y^{\beta}, z^{\gamma}$ are either both out-neighbours or both in-neighbours of 
$x^{\alpha}$. It follows that each vertex of $\tilde{K}_{D}$ is either a source or
a sink in the orientation. Let $X$ the set of sources and $Y$ be the set of sinks.
Then $(X,Y)$ is a bipartition of $\tilde{K}_{D}$. The orientation of $\tilde{K}_{D}$ is from $X$ to $Y$, which transforms into the orientation of $U(D)$ in such a way 
that an edge $xy$ is oriented from $x$ to $y$ if $x \prec y$. Therefore the 
transformed orientation of $U(D)$ is acyclic. 
  
Conversely, assume that there exists a bipartition $(X,Y)$ of $\tilde{K}_{D}$ such 
that the orientation of $\tilde{K}_{D}$ from $X$ to $Y$ transforms into an acyclic 
orientation of $U(D)$. Let $\prec$ be a vertex ordering of $D$ such that $xy$ is
an oriented edge in the orientation of $U(D)$ implies $x \prec y$. 
We prove that $\prec$ is a comparability ordering of $D$. 
Let $x \prec y \prec z$ be any three vertices. Suppose that
$xy$ and $yz$ are arcs of $D$ but $xz$ is not. Hence $(y,x)$ and $(y,z)$ 
are knotted pairs. These two pairs correspond to edges $y^{\beta}x^{\alpha}$ and
$y^{\beta}z^{\gamma}$ in $\tilde{K}_{D}$ where $x^{\alpha}, y^{\beta}, z^{\gamma}$ 
are partial copies of $x, y, z$ respectively. If $y^{\beta} \in X$ then 
the edges $y^{\beta}x^{\alpha}, y^{\beta}z^{\gamma}$ are both oriented from
$y^{\beta}$. Thus the oriented edges $y^{\beta}x^{\alpha}, y^{\beta}z^{\gamma}$
transform into the oriented edges $yx, yz$, which implies $y \prec x$ and 
$y \prec z$, a contradiction. On the other hand, if $y^{\beta} \in Y$, then the 
edges $y^{\beta}x^{\alpha}, y^{\beta}z^{\gamma}$ are both oriented to $y^{\beta}$, 
which transform into the oriented edges $xy, zy$. We then have 
$x \prec y$ and $z \prec y$, a contradiction. So $xz$ must be an arc of $D$. 
A similar proof shows if $zy$ and $yx$ are arcs of $D$ then $zx$ is also an arc of 
$D$. Therefore $\prec$ is a comparability ordering of $D$ and $D$ is 
a comparability digraph.
\qed

\section{The Triangle Lemma for semicomplete digraphs}
\label{semicomplete}

Semicomplete comparability digraphs generalize comparability graphs (see Section~1). 
One of the important tools for the device of an efficient recognition algorithm for
comparability graphs is the {\em Triangle Lemma} (cf. \cite{golumbic}), which is 
a property about implication classes of graphs. In this section, we extend the 
Triangle Lemma for graphs to semicomplete digraphs. We will apply the Triangle Lemma
for semicomplete digraphs to prove that if an implication class of a semicomplete 
digraph does not contain a circuit of length 2, then it contains no circuit at all. 
The Triangle Lemma for semicomplete digraphs will be used in Section~\ref{recog} 
for the device of an efficient recognition algorithm for semicomplete comparability 
digraphs. 

We begin with a simple but useful property of semicomplete digraphs.

\begin{lem}\label{pro:GamDirTri}
Let $D=(V,A)$ be a semicomplete digraph with three distinct vertices $x, y, z$. Then
$(x, y)\Gamma (x, z)$ if and only if the arc between $y$ and $z$ is non-symmetric
and the subdigraph induced by $\{x, y, z\}$ contains a directed triangle.
\end{lem}
\pf Suppose that $(x, y)\Gamma (x, z)$. Then either $yx, xz \in A$ and 
$yz \notin A$ or $zx, xy \in A$ and $zy \notin A$. In the former case 
$zy$ is a non-symmetric arc and $yxz$ is a directed triangle and in the latter case
$yz$ is a non-symmetric arc and $zxy$ is a directed triangle. 
Conversely, suppose that the subdigraph induced by $\{x, y, z\}$ contains 
a directed triangle. If the arc between $y, z$ is non-symmetric then 
either $yx, xz \in A$ and $yz \notin A$ or $zx, xy \in A$ and $zy \notin A$.
In either case, we have $(x, y)\Gamma (x, z)$. 
\qed

The following lemma will be refered to as the {\em Triangle Lemma} for semicomplete 
digraphs due to its resemblance to the Triangle Lemma for graphs, 
cf. \cite{golumbic}. 

\begin{lem}\label{lem:Triangle} (The Triangle Lemma) 
Let $D=(V,A)$ be a semicomplete digraph with three distinct vertices $x, y, z$. 
Let $X, Y, Z$ be the implication classes of $D$ that contain $(y,z), (x,z), (x,y)$
respectively. Suppose that $X, Z$ are nontrivial and $X\neq Z, Z^{-1}, Y$. Then 
for any $(y',z') \in X$, the following two statements hold: 
\begin{itemize}
\item $(x,y')\in Z$ and $(x,z')\in Y$;
\item the arcs between $x$ and $y, y', z, z'$ have the same form, that is, 
      if there is an arc from $x$ to any one of them then there is an arc
      from $x$ to each of them, and if there is an arc from any one of them to $x$,
      then there is an arc from each of them to $x$.
\end{itemize}
\end{lem}
\pf We prove that the two statements hold for each pair $(y',z')$ with 
$(y',z') \neq (y,z)$ and $(y',z') \Gamma (y,z)$. For an arbitrary pair 
$(y',z') \in X$ the validity of the two statements follows easily by induction on 
the length of a $\Gamma$-chain connecting $(y,z)$ and $(y',z')$. So assume that 
$(y',z')$ is such a pair. 
We first show that the arcs between $x$ and $y, z$ have the same form. 
Suppose this is not the case. Then either $zx, xy \in A$ and at least one of $zx, xy$
is non-symmetric or $yx,xz \in A$ and at least one of $yx,xz$ is non-symmetric.
Since these two cases can be transformed to each other by reversing the arcs of $D$,
we only consider the former, i.e., $zx, xy \in A$ and at least one of $zx, xy$ is 
non-symmetric. If $yz\in A$, then $xyz$ is a directed triangle. Since at least one 
of $zx, xy$ is non-symmetric, by Lemma \ref{pro:GamDirTri}, 
$(y,z) \Gamma (y,x)$ or $(y,z) \Gamma (x,z)$, which implies $X = Z^{-1}$ or $X = Y$,
contradicting the assumptions $X \neq Z^{-1}$ and $X \neq Y$. Hence 
$yz \notin A$ and thus $zy$ is a non-symmetric arc.
 
Since $(y',z') \neq (y,z)$ and $(y',z') \Gamma (y,z)$, either $y' \neq y$ and 
$z' = z$ or  $y' = y$ and $z' \neq z$. Suppose first that $y' \neq y$ and $z' = z$.
By Lemma \ref{pro:GamDirTri}, $y'zy$ is a directed triangle and $yy'$ is a 
non-symmetric arc. If $xy' \notin A$ then $y'x$ is a non-symmetric arc. We have 
$(y',z)\Gamma (y',y)\Gamma (x,y)$, which implies $X = Z$, a contradiction to
the assumption $X \neq Z$. So $xy' \in A$. If $xy'$ is a non-symmetric arc, then
$(y',z)\Gamma (x,z)$, which implies $X = Y$, a contradiction to the assumption 
$X \neq Y$. So $xy'$ is a symmetric arc. If $zx$ is non-symmetric, then
$(y',z)\Gamma (y',x)\Gamma (y,x)$, which implies $X = Z^{-1}$, a contradiction to
the assumption $X \neq Z^{-1}$. So $zx$ is a symmetric arc. Since at least one of 
$zx, xy$ is non-symmetric, $xy$ is non-symmetric. We have
$(y',z)\Gamma (y',y)\Gamma (y',x)\Gamma (y,x)$, which implies $X = Z^{-1}$, 
a contradiction. 

Suppose next that $y' = y$ and $z' \neq z$. By Lemma \ref{pro:GamDirTri}, $zyz'$ is 
a directed triangle and $z'z$ is a non-symmetric arc. If $xz' \notin A$ then $z'x$
is a non-symmetric arc. We have $(y,z')\Gamma (y,x)$ which implies $X = Z^{-1}$,
a contradiction. So $xz' \in A$. If $zx$ is a non-symmetric arc, then 
$(y,z')\Gamma (z,z')\Gamma (x,z')\Gamma (x,z)$, which implies $X = Y$, 
a contradiction. So $zx$ is a symmetric arc and thus $xy$ is a non-symmetric arc.   
The arc $xz'$ must be non-symmetric as otherwise 
$(y,z')\Gamma (x,z')\Gamma (x,z)$, which again implies $X = Y$, a contradiction. 
Since $Z$ is nontrivial, there exists a pair $(u,v)$ such that $(u,v) \neq (x,y)$
and $(u,v)\Gamma (x,y)$. Either $u = x$ and $v \neq y$ or $u \neq x$ and $v = y$.
Assume that $u = x$ and $v \neq y$. By Lemma \ref{pro:GamDirTri}, 
$vxy$ is a directed triangle and $yv$ is a non-symmetric arc. 
Note that $v \notin \{x,z,z'\}$. If $vz \in A$ then 
$(y,z)\Gamma (v,z)\Gamma (v,y)\Gamma (v,x)$, which implies $X = Z^{-1}$, 
a contradiction. So $zv$ is a non-symmetric arc. But then 
$(y,z')\Gamma (z,z') \Gamma (z,x)\Gamma (v,x)$, which implies $X = Z^{-1}$, 
a contradiction. Assume on the other hand that $u \neq x$ and $v = y$. 
By Lemma \ref{pro:GamDirTri}, $xyu$ is a directed triangle and $ux$ is 
a non-symmetric arc. We must have $zu \in A$ as otherwise $uz$ is a non-symmetric 
arc and $(y,z)\Gamma (y,u)\Gamma (y,x)$ which implies $X = Z^{-1}$,
a contradiction. If $zu$ is non-symmetric, then 
$(y,z')\Gamma (z,z')\Gamma (z,x)\Gamma (u,x)\Gamma (u,y)\Gamma (x,y)$, which 
implies $X = Z$, a contradiction. So $zu$ is a symmetric arc. But then
$(y,z')\Gamma (z,z')\Gamma (z,x)\Gamma (z,u)\Gamma (y,u)\Gamma (y,x)$, which again 
implies $X = Z^{-1}$, a contradiction.
Therefore the arcs between $x$ and $y,z$ have the same form.
This fact will be used multiple times in the rest of the proof (for different sets 
of three vertices that satisfy similar assumptions as $x, y, z$). 

We now show that the arcs between $x$ and $y',z'$ have the same form as the arcs
between $x$ and $y,z$ and that $(x,y') \in Z$ and $(x,z') \in Y$. We will consider 
three cases depending on the form of the arcs between $x$ and $y,z$. Consider first 
the case when the arcs between $x$ and $y,z$ are all symmetric. 
As the above, either $y' = y$ and $z' \neq z$ or $y' \neq y$ and $z' = z$.
Assume that $y' = y$ and $z' \neq z$. By Lemma \ref{pro:GamDirTri}, 
either $z'yz$ is a directed triangle and $zz'$ is a non-symmetric arc or $zyz'$ is 
a directed triangle and $z'z$ is a non-symmetric arc.  
Assume that $z'yz$ is a directed triangle and $zz'$ is a non-symmetric arc.
(A similar proof applies to the other case.) We must have $z'x \in A$ as otherwise
$xz'$ is a non-symmetric arc and $(y',z')\Gamma (y,x)$, which implies $X = Z^{-1}$, 
a contradiction. Consider the three vertices $x, y', z'$. Clearly, $X$ contains
$(y,z')$, $Z$ contains $(x,y)$, and $Y$ contains $(x,z')$ as $(x,z')\Gamma (x,z)$.
So the three vertices $x, y', z'$ can play the roles of $x, y, z$. It follows from
the above proof that the arcs between $x$ and $y',z'$ must have the same form. 
Since the arcs between $x$ and $y$ are symmetric, the arcs between $x$ and $z'$ are 
also symmetric. Hence the arcs between $x$ and $y,z,z'$ are all symmetric and 
furthermore $(x,z') \in Y$. 
Assume now that $y' \neq y$ and $z' = z$. Again by Lemma \ref{pro:GamDirTri} either 
$yzy'$ is a directed triangle and $y'y$ is a non-symmetric arc or $y'zy$ is 
a directed triangle and $yy'$ is a non-symmetric arc. Assume that 
$yzy'$ is a directed triangle and $y'y$ is a non-symmetric arc. (The proof for 
the other case is similar.) We must have $xy' \in A$ as otherwise 
$y'x$ is a non-symmetric arc and $(y',z)\Gamma (x,z)$, which implies $X = Y$, 
a contradiction. Thus $(x,y')\Gamma (x,y)$ and so $(x,y') \in Z$. Now the vertices 
$x, y', z$ can play the roles of $x, y, z$. We know from above that the arcs between
$x$ and $y',z$ have the same form. Since the arcs between $x$ and $z$ are symmetric,
the arcs between $x$ and $y'$ are also symmetric. Therefore the arcs between $x$ 
and $y,y',z$ are all symmetric and furthermore, $(x,y') \in Z$.  

Suppose next that the arcs between $x$ and $y,z$ are non-symmetric and from $x$ to 
$y,z$. We only consider the case when $y' = y$, $z'yz$ is a directed triangle and 
$zz'$ is a non-symmetric arc. (Proofs for other cases are similar.)
If $z'x \in A$ then $(x,z')\Gamma (x,z)$ and so $(x,z')\in Y$. The three vertices 
$x,y,z'$ can now play the roles of $x, y, z$. We know from above that the arcs 
between $x$ and $y,z'$ have the same form. This contradicts the fact that the arcs 
between $x$ and $y,z'$ do not have the same form. So $z'x \notin A$, that is, 
$xz'$ is a non-symmetric arc. Hence the arcs between $x$ and $y,z,z'$ are all
non-symmetric and from $x$ to $y,z,z'$. It remains to show that $(x,z') \in Y$. 
Since $Z$ is nontrivial, there exists a pair $(u,v)$ with $(u,v) \neq (x,y)$
and $(u,v)\Gamma (x,y)$. Either $u = x$ and $v \neq y$ or $u \neq x$ and $v = y$.
Assume that $u = x$ and $v \neq y$. (A similar proof applies when $u \neq x$ and 
$v = y$.) By Lemma \ref{pro:GamDirTri}, $vxy$ is a directed triangle and $yv$ is 
a non-symmetric arc.
If $vz'$ is a non-symmetric arc, then $(y,z')\Gamma (y,v)\Gamma (x,v)$, which
implies $X = Z$, a contradiction. If $vz'$ is a symmetric arc, then 
$(y,z')\Gamma (v,z')\Gamma (v,x)$ which implies $X = Z^{-1}$, a contradiction. 
So $z'v$ is a non-symmetric arc. Assume $vz \in A$. If $vz$ and $yz$ are both 
symmetric, then $(y,z)\Gamma (v,z)\Gamma (v,x)$, which implies $X = Z^{-1}$, 
a contradiction.  If $vz$ and $yz$ are both non-symmetric, then
$(y,z')\Gamma (z,z')\Gamma (v,z')\Gamma (v,x)$, which implies $X = Z^{-1}$,
a contradiction. If $vz$ is symmetric and $yz$ is non-symmetric, then
$(y,z')\Gamma (z,z')\Gamma (z,v)\Gamma (x,v)$ which implies $X = Z$,
a contradiction. Finally if $vz$ is non-symmetric and $yz$ is symmetric, then
$(y,z)\Gamma (y,v)\Gamma (x,v)$ which implies $X = Z$, a contradiction.
Hence $vz \notin A$ and $zv$ is a non-symmetric arc. Therefore
$(x,z')\Gamma (x,v)\Gamma (x,z)$ which means $(x,z') \in Y$ and we are done.

Suppose now that the arcs between $x$ and $y,z$ are non-symmetric and from 
$y,z$ to $x$. Once again we only consider the case when $y' = y$, $z'yz$ is 
a directed triangle and $zz'$ is a non-symmetric arc. Suppose first that $yz$ is 
a non-symmetric arc. We must have $z'x \in A$ as otherwise $xz'$ is a non-symmetric
arc and $(y,z')\Gamma (y,x)$, which implies $X = Z^{-1}$, a contradiction. 
Since $Z$ is nontrivial, there exists a pair $(u,v)$ with $(u,v) \neq (x,y)$ 
and $(u,v)\Gamma (x,y)$. Consider first the case when $u = x$,
$yxv$ is a directed triangle and $vy$ is a non-symmetric arc. 
If $zv \in A$, then $(y,z)\Gamma (v,z)\Gamma (v,y)\Gamma (v,x)$, which implies 
$X = Z^{-1}$, a contradiction. So $zv \notin A$ and $vz$ is a non-symmetric arc. 
If $z'v \in A$ then $(y,z')\Gamma (z,z')\Gamma (v,z')\Gamma (v,z)\Gamma (v,x)$,
which implies $X = Z^{-1}$, a contradiction. So $z'v \notin A$ and $vz'$ is a
non-symmetric arc. We have $(x,z')\Gamma (x,v)\Gamma (x,z)$ and so 
$(x,z') \in Y$. Thus the three vertices $x,y,z'$ can play the roles of $x,y,z$.
We know from above that the arcs between $x$ and $y,z'$ have the same form. 
Since $yx$ is a non-symmetric arc, $z'x$ is a non-symmetric arc. Therefore the arcs
between $x$ and $y,z,z'$ have the same form. Consider now the case when $v = y$, 
$uyx$ is a directed triangle and $xu$ is a non-symmetric arc. 
We must have $uz \in A$ as otherwise $zu$ is a non-symmetric arc and 
$(y,z)\Gamma (y,u)$, which implies $X = Z^{-1}$, a contradiction. Similarly,
we must have $uz' \in A$ as otherwise $z'u$ is a non-symmetric arc and 
$(y,z')\Gamma (z,z')\Gamma (z,u)\Gamma (x,u)\Gamma (y,u)$, which implies 
$X = Z^{-1}$, a contradiction. We know from above that $z'x \in A$. If $z'x$ is
symmetric, then $uz'$ is also symmetric as otherwise 
$(y,z')\Gamma (x,z')\Gamma (x,u)\Gamma (y,u)$, which implies
$X = Z^{-1}$, a contradiction. The vertices $u,x,z'$ can play the roles of 
$x,y,z$ respectively. Indeed, $(x,z')\Gamma (y,z') \in X$,    
$(u,z')\Gamma (u,z)\Gamma (x,z) \in Y$, and $(u,x)\Gamma (u,y) \in Z$. So
the arcs between $u$ and $x,z'$ have the same form. But this contradicts the fact
that $xu$ is a non-symmetric arc and $uz'$ is symmetric. Hence $z'x$ is a 
non-symmetric arc and therefore the arcs between $x$ and $y,z,z'$ have the same form.
Moreover, $(x,z')\Gamma (u,z')\Gamma (u,x)\Gamma (u,z)\Gamma (x,z)$, which implies
$(x,z') \in Y$. 

Suppose then that the arc $yz$ is symmetric. We must have $z'x \in A$ as otherwise
$xz'$ is a non-symmetric arc and $(y,z')\Gamma (y,x)$, which implies $X = Z^{-1}$,
a contradiction. Since $Z$ is nontrivial, there exists a pair $(u,v)$ with 
$(u,v) \neq (x,y)$ and $(u,v)\Gamma (x,y)$. Consider first the case when $u = x$,
$yxv$ is a directed triangle and $vy$ is a non-symmetric arc.
If $vz \notin A$ then $zv$ is a non-symmetric arc and 
$(y,z)\Gamma (y,v)\Gamma (x,v)$ which implies $X = Z$, a contradiction. So
$vz \in A$. The arc $vz$ must be non-symmetric as otherwise 
$(y,z)\Gamma (v,z)\Gamma (v,x)$ which implies $X = Z^{-1}$, a contradiction. 
Suppose that $xz' \in A$ (i.e., it is symmetric).
Then $z'y$ is symmetric as otherwise $(y,z')\Gamma (x,z')\Gamma (x,y)$ which
implies $X = Z$, a contradiction. If $z'v$ is a non-symmetric arc then 
$(y,z')\Gamma (y,v)\Gamma (x,v)$ which implies $X = Z$, a contradiction. 
If $z'v$ is a symmetric arc then $(y,z')\Gamma (v,z')\Gamma (v,z)\Gamma (v,x)$
which implies $X = Z^{-1}$, a contradiction. So $vz'$ is a non-symmetric arc. 
Thus $(x,z')\Gamma (x,v)\Gamma (x,z)$ and in particular $(x,z') \in Y$.
We know from above that the arcs between $x$ and $y,z'$ have the same form. 
This contradicts the fact that the arc between $x$ and $y$ is non-symmetric and 
the arcs between $x$ and $z'$ are symmetric. Hence $xz' \notin A$ and $z'x$ is 
a non-symmetric arc. Therefore the arcs between $x$ and $y,z,z'$ have the same form. 
Since $(x,z)\Gamma (x,v)$, $Y = Z$. Since $X \neq Z^{-1}$, $X \neq Y^{-1}$.
If $z'v \in A$, then $(v,z')\Gamma (v,z)\Gamma (v,x)\Gamma (z,x) \in Y^{-1}$. 
This is a contradiction because with vertices $v, y, z'$ we have
$(y,z') \in X$, $(v,z') \in Y^{-1}$ and $(v,y) \in Z^{-1}$. From the above 
we know that the arcs between $v$ and $y,z'$ have the same form. But this is not 
the case: $vy$ is a non-symmetric arc and $z'v \in A$. So $vz'$ must be
a non-symmetric arc and hence $(x,z') \Gamma (x,v)\Gamma (x,z) \in Y$. 

Consider now the case when $v = y$, $uyx$ is a directed triangle and $xu$ is 
a non-symmetric arc. As above we have $z'x \in A$. We must also have $uz \in A$ 
as otherwise $zu$ is a non-symmetric arc and $(y,z)\Gamma (y,u)$, which implies 
$X = Z^{-1}$, a contradiction. Thus $(u,z)\Gamma (x,z)$ which means $(u,z) \in Y$.
The three vertices $u,y,z$ can play the roles of $x,y,z$. Hence the arcs between 
$u$ and $y,z$ have the same form. Either the arcs between $u$ and $y,z$ are
symmetric or non-symmetric and from $u$ to $y,z$. Suppose that the arcs between $u$ 
and $y,z$ are non-symmetric and from $u$ to $y,z$. If $z'u \in A$ then 
$(u,z')\Gamma (u,z)\Gamma (x,z)$. In particular, $(u,z') \in Y$. So
the arcs between $u$ and $y,z'$ have the same form. But this is not true as
$yu \notin A$ and $z'u \in A$. So $uz'$ is a non-symmetric arc.
Since $(x,z')\Gamma (x,u)\Gamma (x,z)$, $(x,z') \in Y$. Now $x,y,z'$ can play 
the roles of $x,y,z$. The arcs between $x$ and $y,z'$ have the same form. 
Since $yx$ is non-symmetric, $z'x$ is also non-symmetric.
Therefore the arcs between $x$ and $y,z,z'$ have the same form.
Suppose that the arcs between $u$ and $y,z$ are symmetric. We must have $z'u \in A$
as otherwise $uz'$ is a non-symmetric arc and $(y,z')\Gamma (y,u)$, which implies 
$X = Z^{-1}$, a contradiction. Note that $(u,z')\Gamma (u,z)\Gamma (x,z) \in Y$.
Now $u,y,z'$ can play the roles of $x,y,z$. We conclude that the arcs between $u$ 
and $y,z'$ have the same form. Since the arcs between $u$ and $y$ are symmetric, 
the arcs between $u$ and $z'$ are also symmetric. 
If $z'x$ is symmetric then $(y,z')\Gamma (x,z')\Gamma (u,z')\Gamma (u,z)\Gamma
(x,z)$ which implies $X = Y$, a contradiction. So $z'x$ is a non-symmetric arc. 
Therefore the arcs between $x$ and $y,z,z'$ have the same form. Furthermore 
$(x,z')\Gamma (u,z')\Gamma (u,z)\Gamma (x,z)$ and so $(x,z') \in Y$. 
\qed

For convenience we shall use $[x,y]$ to denote the implication class containing
the pair $(x,y)$. 

\begin{lem}\label{lem:path-chord}
Let $D$ be a semicomplete digraph and $I\in \mathcal{I}(D)$ with $I\neq I^{-1}$.
Suppose that $x_{0}, x_{1}, \dots, x_{k-1}$ are distinct vertices such that 
$(x_i,x_{i+1})\in I$ for each $0 \leq i < k-1$ and $(x_i,x_j) \notin I$ for all
$i,j$ with $i-j\neq \pm1$. Then the following statements hold:
\begin{itemize}
\item For all $i,j$ with $i-j\neq \pm1$, $[x_i,x_j]$ is a trivial implication class.
\item If $x_ix_{i+1}\in A$ for each $0 \leq i < k-1$, then $x_ix_j\in A$ for all 
      $i<j$.
\item If $x_{i+1}x_i\in A$ for each $0 \leq i < k-1$, then $x_ix_j\in A$ for all 
      $i>j$.
\end{itemize}
\end{lem}
\pf For the first statement assume to the contrary that $[x_i,x_j]$ is a nontrivial 
implication class for some $i,j$ with $i-j \neq \pm1$. Choose such $i,j$ so that 
$i<j$ and $j-i$ is the smallest. We apply the Triangle Lemma with
$x,y,z$ corresponding to $x_i,x_j,x_{j-1}$ respectively. By assumption 
$[x_j,x_{j-1}] = I^{-1}$, $[x_i,x_j] \neq I$ and $[x_i,x_j] \neq I^{-1}$. 
When $j=i+2$, $[x_i,x_{j-1}] = [x_i,x_{i+1}] = I \neq I^{-1}$. When $j > i+2$,
the choice of $i,j$ implies $[x_i,x_{j-1}]$ is trivial so 
$[x_i,x_{j-1}] \neq I^{-1}$. Since $(x_{j-1},x_{j-2}) \in I^{-1}$, 
by the Triangle Lemma $(x_i,x_{j-1}) \in [x_i,x_j]$. This contradicts the choice of
$i,j$.

The proofs for the second and the third statements are similar so we only prove the 
second one. Suppose that $x_ix_j \notin A$ for some $i<j$. Choose such $i, j$ so 
that $i<j$ and $j-i$ is the smallest. The choice of $i,j$ implies that $j-i \geq 2$,
$x_jx_i$ is a non-symmetric arc, and $x_ix_{j-1}x_j$ is a directed triangle.
By Lemma \ref{pro:GamDirTri} $(x_{j-1},x_i)\Gamma (x_{j-1},x_j)$ and so
$[x_{j-1},x_i] = [x_{j-1},x_j] = I$. If $j-i=2$, then
$I = [x_{j-1},x_i] = [x_{i+1},x_i] = I^{-1}$, which contradicts the assumption that
$I\neq I^{-1}$. If $j-i>2$, then by the first statement $[x_{j-1},x_i]$ is trivial,
a contradiction to the fact that $I$ is nontrivial.
\qed

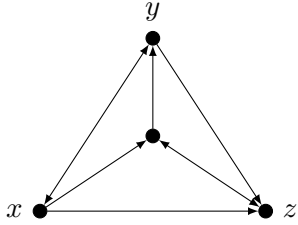
\begin{figure}[ht]
\begin{center}
\begin{tikzpicture}[>=latex]
                \node [label={above:$y$}] [style=blackvertex] (1) at (0,2.3) {};
                \node [style=blackvertex] (2) at (0,1) {};
                \node [label={left:$x$}] [style=blackvertex] (3) at (-1.5,0) {};
                \node [label={right:$z$}] [style=blackvertex] (4) at (1.5,0) {};

                \draw [-> ] (3) to (2);
                \draw [-> ] (2) to (1);
                \draw [-> ] (1) to (4);
                \draw [-> ] (3) to (4);
                \draw [<->] (1) to (3);
                \draw [<->] (2) to (4);
\end{tikzpicture}
\end{center}
\vspace{-2mm}
\caption{A semicomplete comparability digraph $D$ for which not all $I\in \mathcal{I}(D)$ are transitive.}
\label{non-trans}
\end{figure}

Using the Triangle Lemma for graphs one can prove easily that any implication class 
of a graph that contains no circuit of length 2 is transitive and hence contains no 
circuit at all, cf. \cite{golumbic}. This property is however not true for digraphs.
The semicomplete digraph in Figure \ref{non-trans} is a comparability digraph so
no implication class contains a circuit but $[x,y]$ is not transitive 
(as it contains $(x,y), (y,z)$ but not $(x,z)$). In general, we have the following, 
which is the case $k=2$ of Lemma \ref{lem:path-chord}.

\begin{cor} \label{lem:IIInvTri}
Let $D$ be a semicomplete digraph and $I\in \mathcal{I}(D)$ with $I\neq I^{-1}$.
Then for any three vertices $x, y, z$ of $D$, if $(x,y), (y,z)\in I$ then
$(x,z) \in I$, or $(z,x) \in I$, or $[x,z]$ is trivial.
\qed
\end{cor}

Let $D$ be a digraph and $I \in \mathcal{I}(D)$. Let 
$(x_{0},x_{1}), (x_{1},x_{2}), \dots, (x_{k-2},x_{k-1}), (x_{k-1},x_0)$ be 
a circuit in $I$. The circuit is called {\em induced} in $I$ if $(x_i,x_j) \notin I$
for all $i, j$ with $|i-j| \notin \{1,k-1\}$. A pair $(x_i,x_j)$ with 
$|i-j| \notin \{1,k-1\}$ is called a {\em chord} of the circuit. Note that 
the vertices of an induced circuit are distinct. If $(x_i,x_j)$ is a chord then 
$(x_j,x_i)$ is also a chord. A chord exists only when $k \geq 4$. 

\begin{cor}\label{cor:circuit-chord}
Let $D$ be a semicomplete digraph and $I\in \mathcal{I}(D)$ with $I\neq I^{-1}$.
Suppose that $(x_{0},x_{1}), (x_{1},x_{2}), \dots, (x_{k-2},x_{k-1}), (x_{k-1},x_0)$
is an induced circuit in $I$. Then for any chord $(x_i,x_j)$ of the circuit, 
$[x_i,x_j]$ is a trivial implication class.
\end{cor}
\pf Since the circuit is induced in $I$ and $(x_i,x_j)$ is a chord, 
$(x_i,x_j) \notin I$ and $(x_i,x_j) \notin I^{-1}$.
By Lemma \ref{lem:path-chord} on the vertices $x_i, x_{i+1}, \dots, x_j$, 
$[x_i,x_j]$ is a trivial implication class. 
\qed

\begin{lem}\label{lem:I-no-circuit}
Let $D = (V,A)$ be a semicomplete digraph and $I\in \mathcal{I}(D)$ with 
$I\neq I^{-1}$. Assume that $C:\ (x_{0},x_{1}), (x_{1},x_{2}), \dots, 
(x_{k-2},x_{k-1}), (x_{k-1},x_0)$ is an induced circuit in $I$. Then the subdigraph 
of $D$ induced by $x_0, x_1, \dots, x_{k-1}$ is a complete digraph.  
\end{lem}
\pf Since $I \neq I^{-1}$, $k \geq 3$. Suppose that $x_0x_1\dots x_{k-1}$ is a 
directed cycle. When $k=3$, $x_0x_1x_2$ is a directed triangle. It follows from 
Lemma \ref{pro:GamDirTri} and the assumption $I \neq I^{-1}$ that every arc in the 
directed triangle is symmetric. So the subgraph of $D$ induced by $x_0,x_1,x_2$ is 
a complete digraph. Assume $k > 3$. Let $(x_i,x_j)$ be a chord of the circuit $C$. 
By Corollary \ref{cor:circuit-chord}, $[x_i,x_j]$ is a trivial implication class.
Since $x_0x_1\dots x_{k-1}$ is a directed cycle, Lemma \ref{lem:path-chord} on 
$x_i,x_{i+1},\dots,x_j$ and $x_j,x_{j+1},\dots,x_i$ implies that 
$x_ix_j, x_jx_i \in A$, i.e., the arcs between $x_i$ and $x_j$ are symmetric. 
If any $x_ix_{i+1}$ is non-symmetric, then $(x_i,x_{i+2})\Gamma (x_{i+1},x_{i+2})$, 
which implies $(x_i,x_{i+2}) \in I$, contradicting that $C$ is an induced circuit
in $I$. So the arcs between $x_i$ and $x_{i+1}$ are symmetric for
all $i$. Hence the subdigraph of $D$ induced by $x_0,x_1,\dots,x_{k-1}$ is
a complete digraph. 

It remains to show that $x_0x_1 \dots x_{k-1}$ is a directed cycle. Suppose not. 
Then there are two non-symmetric arcs in opposite directions along the circuit $C$. 
By relabeling the vertices and/or swapping $I$ and $I^{-1}$ if necessary, we may 
assume that $x_1x_0$ and $x_tx_{t+1}$ are non-symmetric arcs for some 
$1\leq t\leq k-1$. Moreover, we may assume that, for each $1\leq i\leq t-1$, 
the arcs between $x_i$ and $x_{i+1}$ are symmetric.  

Consider first the case when $t = 1$ (i.e., $x_1x_2$ is a non-symmetric arc). 
Since $I$ is nontrivial and $(x_0,x_1)\in I$, there exists a vertex $u$ such that 
$(x_0,x_1)\Gamma (x_0,u)$ or $(x_0,x_1)\Gamma (u,x_1)$. 
By Lemma \ref{pro:GamDirTri}, $x_0ux_1$ is a directred cycle and one of $x_0u, ux_1$
is a non-symmetric arc. If $x_0u, ux_1$ are both non-symmetric, then 
$(x_0,x_1)\Gamma (x_0,u)\Gamma (x_1,u)\Gamma (x_1,x_0)$, which implies $I = I^{-1}$,
a contradiction to the assumption $I \neq I^{-1}$. So one of $x_0u, ux_1$ is 
non-symmetric and the other is symmetric. Assume that $ux_1$ is non-symmetric; thus 
$x_0u$ is symmetric. If $x_2u \in A$ then 
$(x_0,x_1)\Gamma (x_0,u)\Gamma (x_1,u)\Gamma (x_2,u)\Gamma (x_2,x_1)$, which implies
$I = I^{-1}$, a contradiction. So $ux_2$ is a non-symmetric arc.
If $x_2x_0 \in A$ then $(x_0,x_1)\Gamma (x_0,u)\Gamma (x_0,x_2)$, which implies 
$(x_0,x_2) \in I$. This contradicts that $I \neq I^{-1}$ and $C$ is an induced 
circuit in $I$. So $x_0x_2$ is a non-symmetric arc. Similarly, there exists 
a vertex $v$ such that $x_1x_2v$ is a directed cycle and exactly one of $x_2v, vx_1$
is non-symmetric. Assume that $x_2v$ is non-symmetric; thus $vx_1$ is symmetric.
If $vx_0 \in A$, then $(x_1,x_2)\Gamma (x_1,v)\Gamma (x_2,v)\Gamma (x_0,v) \Gamma
(x_0,x_2)$, which again implies $(x_0,x_2) \in I$, a contradiction. So $x_0v$ is 
a non-symmetric arc. But then $(x_0,x_1)\Gamma (v,x_1)\Gamma (x_2,x_1)$, which
implies $I = I^{-1}$, a contradiction. Assume that $vx_1$ is non-symmetric; thus
$x_2v$ is symmetric. If $x_0v \in A$, then 
$(x_1,x_2)\Gamma (v,x_2)\Gamma (v,x_1)\Gamma (v,x_0)\Gamma (x_1,x_0)$,
which implies $I = I^{-1}$, a contradiction. So $vx_0$ is a
non-symmetric arc. But then $(x_1,x_2)\Gamma (v,x_2)\Gamma (x_0,x_2)$, which
implies $(x_0,x_2)\in I$, a contradiction.

Assume that $x_0u$ is non-symmetric; thus $ux_1$ is symmetric. If $x_2u$ is 
a non-symmetric arc, then $(x_0,x_1)\Gamma (u,x_1)\Gamma (x_2,x_1)$, which implies
$I = I^{-1}$, a contradiction. So $ux_2 \in A$. If $x_2x_0$ is a non-symmetric arc,
then $(x_0,x_1)\Gamma (u,x_1)\Gamma (u,x_0)\Gamma (u,x_2)\Gamma (x_0,x_2)$, which
implies $(x_0,x_2) \in I$. This contradicts that $I \neq I^{-1}$ and $C$ is 
an induced circuit in $I$. So $x_0x_2 \in A$. As above, there exists
a vertex $v$ such that $x_1x_2v$ is a directed cycle and exactly one of $x_2v, vx_1$
is non-symmetric. Assume that $vx_1$ is non-symmetric; thus $x_2v$ is symmetric. 
If $x_0v \in A$, then $(x_0,x_1)\Gamma (x_0,v)\Gamma (x_1,v)\Gamma (x_2,v)\Gamma
(x_2,x_1)$, which implies $I = I^{-1}$, a contradiction. So $vx_0$ is a non-symmetric
arc. But then $(x_1,x_2)\Gamma (v,x_2)\Gamma (x_0,x_2)$, which implies 
$(x_0,x_2) \in I$, a contradiction.
Assume now that $x_2v$ is non-symmetric; thus $vx_1$ is symmetric.
If the arc between $u$ and $v$ is non-symmetric, then $(x_0,x_1)\Gamma (u,x_1)\Gamma
(v,x_1)\Gamma (x_2,x_1)$, which implies $I = I^{-1}$, a contradiction. So 
the arcs between $u, v$ are symmetric. If $x_0v$ is a non-symmetric arc, then
$(x_0,x_1)\Gamma (v,x_1)\Gamma (x_2,x_1)$, which implies $I=I^{-1}$, a contradiction.
So $vx_0 \in A$. If the arc $ux_2$ is non-symmetric, then 
$(x_1,x_2)\Gamma (x_1,v)\Gamma (x_2,v)\Gamma (u,v)\Gamma (x_0,v)\Gamma (x_0,x_2)$,
which implies $(x_0,x_2) \in I$, a contradiction. So $ux_2$ is symmetric. But then 
$(x_0,x_1)\Gamma (u,x_1)\Gamma (u,x_2)\Gamma (u,v)\Gamma (x_0,v)\Gamma (x_0,x_2)$,
which also implies $(x_0,x_2) \in I$, a contradiction. 

A similar proof applies for the case when $t = k-1$. Thus $1 < t < k-1$. By our
assumption above, the arcs between $x_i$ and $x_{i+1}$ are symmetric for each 
$1\leq i\leq t-1$. In particular, when $t=2$, the arcs between $x_1$ and $x_t$ 
are symmetric. Applying Lemma \ref{lem:path-chord} to the vertices 
$x_0, x_1, \dots, x_t$, we have $x_tx_i \in A$ for each $0\leq i \leq t-2$.
Similarly, applying Lemma \ref{lem:path-chord} to the vertices 
$x_1, x_2, \dots, x_{t+1}$, we have $x_1x_i \in A$ for each $2 \leq i \leq t+1$.
In particular, $x_tx_0, x_1x_{t+1} \in A$ and when $t > 2$, the arcs between 
$x_1$ and $x_t$ are symmetric. So in any case the arcs between $x_1$ and $x_t$
are symmetric. According to Corollary \ref{cor:circuit-chord}, $[x_t,x_0]$ and
$[x_1,x_{t+1}]$ are trivial and when $t > 2$, $[x_1,x_t]$ is also trivial.
Hence $x_tx_0$ is non-symmetric as otherwise $(x_t,x_0)\Gamma (x_t,x_1)$, 
a contradiction to that $[x_t,x_0]$ is trivial. Similarly, $x_1x_{t+1}$ is also
non-symmetric. As above, there exists a vertex $u$ such that $x_0ux_1$ is 
a directed triangle and exactly one of $x_0u, ux_1$ is non-symmetric. Assume that 
$x_0u$ is non-symmetric; thus $ux_1$ is symmetric. If $ux_t \in A$, then 
$(x_t,x_0)\Gamma (x_t,u)$, contradicting that $[x_t,x_0]$ is trivial. So $x_tu$ is 
a non-symmetric arc. 
Then $(x_0,x_1)\Gamma (u,x_1)\Gamma (x_t,x_1)$, which implies $(x_t,x_1) \in I$.  
Since $I$ is non-trivial and $[x_t,x_1]$ is trivial when $t > 2$, we have $t = 2$. 
But when $t=2$, $(x_t,x_1) = (x_2,x_1) \in I^{-1}$ so $(x_t,x_1)$ is in both $I$
and $I^{-1}$, a contradiction to the assumption $I \neq I^{-1}$.
Assume that $ux_1$ is non-symmetric; thus $x_0u$ is symmetric. If $ux_t$ is
a non-symmetric arc, then $(x_t,x_0)\Gamma (u,x_0)$, a contradiction to the fact
that $[x_t,x_0]$ is trivial. So $x_tu \in A$. The arc $x_tu$ must be non-symmetric 
as otherwise $(x_0,x_1)\Gamma (x_0,u)\Gamma (x_t,u)\Gamma (x_t,x_1)$, which implies
$(x_t,x_1)\in I$, a contradiction. If $x_{t+1}u \in A$, then $(x_1,x_{t+1})\Gamma
(u,x_{t+1})$, contradicting the fact that $[x_1,x_{t+1}]$ is trivial. So $ux_{t+1}$ 
is a non-symmetric arc. As above, there exists a vertex $v$ such that $x_{t+1}vx_t$ 
is a directed triangle and exactly one of $x_{t+1}v, vx_t$ is non-symmetric.
Assume that $x_{t+1}v$ is non-symmetric; thus $vx_t$ is symmetric.
We must have $vx_1 \notin A$ as otherwise $(x_1,x_{t+1})\Gamma (x_1,v)$, which 
contradicts $[x_1,x_{t+1}]$ is trivial. So $x_1v$ is a non-symmetric arc. 
Then $(x_t,x_{t+1})\Gamma (x_t,v)\Gamma (x_t,x_1)$, which implies $(x_t,x_1)\in I$, 
also a contradiction. Hence $vx_t$ is non-symmetric; thus $x_{t+1}v$ is symmetric. 
Since $[x_1,x_{t+1}]$ is trivial, $x_1v \in A$. If $x_1v$ is symmetric, then
$(x_t,x_{t+1})\Gamma (v,x_{t+1})\Gamma (v,x_1)\Gamma (x_t,x_1)$, which implies 
$(x_t,x_1)\in I$, a contradiction. So $x_1v$ is non-symmetric. 
If $x_0v \in A$, then $(x_t,x_0)\Gamma (v,x_0)$, contradicting that $[x_t,x_0]$
is trivial. So $vx_0$ is a non-symmetric arc. We claim that $(x_0,x_{t+1}) \in I$.
Indeed, if $x_0x_{t+1} \in A$, then 
$(x_0,x_{t+1})\Gamma (v,x_{t+1})\Gamma (x_t,x_{t+1}) \in I$. On the other hand,
if $x_{t+1}x_0$ is a non-symmetric arc, then 
$(x_0,x_{t+1})\Gamma (x_0,u)\Gamma (x_0,x_1) \in I$. Since $C$ is an induced 
circut in $I$, $(x_0,x_{t+1})$ is not chord of $C$. So we have $t = k-2$.
But then $(x_0,x_{t+1}) = (x_0,x_{k-1}) \in I^{-1}$. Thus $(x_0,x_{t+1})$ is in both
$I$ and $I^{-1}$, a contradiction to the fact $I \neq I^{-1}$. 
\qed

\begin{thm} \label{thm:no-circuit-eq}
Let $D$ be a semicomplete digraph and $I\in \mathcal{I}(D)$. Then $I\neq I^{-1}$ if 
and only if $I$ contains no circuit.
\end{thm}
\pf  The sufficiency follows from Proposition \ref{ob:NoCirInvIntEmp}.
For the necessity suppose that $I\neq I^{-1}$ and $I$ contains a circuit. Let 
$C:\ (x_0,x_1), (x_1,x_2), \dots, (x_{k-2},x_{k-1}),(x_{k-1},x_0)$ be an induced 
circuit in $I$. By Lemma \ref{lem:I-no-circuit}, the subdigraph of $D$ induced by 
$x_0, x_1, \dots, x_{k-1}$ is a complete digraph. Since $(x_0,x_1)\in I$ and 
$I$ is nontrivial, $(x_0,x_1) \Gamma (u,v)$ for some $(u,v) \neq (x_0,x_1)$.

{\bf Claim:} $k = 3$ and $(u,v),(v,x_2),(x_2,u)$ is a circuit in $I$.

{\bf Proof of Claim.} Either $u=x_0$ or $v=x_1$. 
When $u=x_0$, by Lemma \ref{pro:GamDirTri}, there is a non-symmetric arc between 
$v$ and $x_1$ and the subdigraph of $D$ induced by $v, x_0, x_1$ contains 
a directed triangle. Similarly, when $v=x_1$, there is a non-symmetric arc between 
$u$ and $x_0$ and the subdigraph of $D$ induced by $u, x_0, x_1$ contains 
a directed triangle.

Consider the case when $u = x_0$. We show that the arcs between $v$ and $x_0$ are 
symmetric. Suppose that $x_0v$ is a non-symmetric arc. (The proof for the case when 
$vx_0$ is non-symmetric is similar.)
Since the subdigraph of $D$ induced by $v, x_0, x_1$ contains a directed triangle, 
$vx_1$ is an arc. Note from above that $vx_1$ is also non-symmetric. If $x_2v$ is 
a non-symmetric arc then $(x_0,x_1)\Gamma (v,x_1)\Gamma (x_2,x_1)$, which implies 
$(x_2,x_1) \in [x_0,x_1] = I$. Thus $(x_2,x_1)$ is in both $I$ and $I^{-1}$, 
contradicting the assumption $I \neq I^{-1}$. So $vx_2$ is an arc. If $vx_2$ is 
non-symmetric, then $(x_0,x_1)\Gamma (x_0,v)\Gamma (x_0,x_2)$, which implies 
$(x_0,x_2) \in I$. 
On the other hand, if $vx_2$ is symmetric, then $(x_1,x_2)\Gamma (v,x_2)\Gamma 
(x_0,x_2)$, which also implies $(x_0,x_2) \in I$. Since $C$ is an induced circuit 
in $I$, $(x_0,x_2)$ is not a chord of $C$. So $k = 3$ and $(x_2,x_0) \in I$. 
Thus $I$ contains both $(x_0,x_2)$ and $(x_2,x_0)$, a contradiction to the 
assumption $I \neq I^{-1}$. Therefore the arcs between $v$ and $x_0$ are symmetric.
If there is a non-symmetric arc between $v$ and $x_{k-1}$, then 
$(x_0,x_1)\Gamma (x_0,v)\Gamma (x_0,x_{k-1})$, which implies $(x_0,x_{k-1}) \in I$,
a contradiction to that $I \neq I^{-1}$. So the arcs between $v$ and $x_{k-1}$ 
are symmetric. Then $(v,x_{k-1})\Gamma (x_1,x_{k-1})$, showing that $[x_1,x_{k-1}]$
is a non-trivial implication class. So we must have $k = 3$ as otherwise it 
contradicts Corollary \ref{cor:circuit-chord}. Since $(v,x_2)\Gamma (x_1,x_2)$,
$(v,x_2) \in I$. Hence $(x_0,v), (v,x_2), (x_2,x_0)$ is a circuit in $I$. 
When $v = x_1$, a similar proof shows that $k = 3$ and $(u,x_1),(x_1,x_2),(x_2,u)$
is a circuit in $I$.
\qed

Since $(x_0,x_1), (x_1,x_2) \in I$, there exists a $\Gamma$-chain connecting
$(x_0,x_1)$ and $(x_1,x_2)$:
\[(x_0,x_1) = (u_0,v_0)\Gamma (u_1,v_1)\Gamma \cdots \Gamma (u_\ell,v_\ell) =
           (x_1,x_2).\]
Applying the above claim repeatedly to circuits $(u_i,v_i),(v_i,x_2),(x_2,u_i)$,
we have $(x_2,u_i) \in I$ for each $0 \leq i \leq \ell$. In particular,
we have $(x_2,x_1) = (x_2,u_\ell) \in I$, which contradicts the assumption
$I \neq I^{-1}$. 
\qed

\section{Recognition of semicomplete comparability digraphs}
\label{recog}

In this section, we describe an algorithm for finding a comparability ordering
of a semicomplete digraph whenever possible. A necessary condition for a digraph $D$
to have a comparability ordering is that $I \neq I^{-1}$ for all 
$I \in \mathcal{I}(D)$ according to Proposition~\ref{pro:NoCirComImp}. 
Our algorithm will indeed succeed for each semicomplete digraph that satisfies this 
necessary condition. Hence the condition is both necessary and sufficient for 
a semicomplete digraph to be a comparability digraph. This is in sharp contrast to 
the fact that there are non-comparability digraphs for which none of the implication 
classes contains a circuit (see example in Figure~\ref{noCs}).

Suppose that $D$ is a semicomplete digraph such that $I \neq I^{-1}$ for all 
$I \in \mathcal{I}(D)$. Our algorithm will construct a set $T \subset Z_D$ by adding
to $T$ all pairs in one of $I, I^{-1}$ for each non-trivial implication class
$I\in \mathcal{I}(D)$. In doing so it gives the priority to those which contain 
pairs implied by the transitivity of $T$. More specifically, it proceeds as follows. 
Initially, $T = \emptyset$ and ${\cal I} = \{I\in \mathcal{I}(D):\
|I| > 1\}$. If there are pairs $(x,y), (y,z) \in T$ such that $(x,z) \in I$ 
for some $I \in \cal I$, then add all pairs of $I$ to $T$ and discard $I$ and 
$I^{-1}$ from $\cal I$; otherwise arbitrarily pick an $I \in \cal I$ and add all 
pairs of $I$ to $T$ and discard $I$ and $I^{-1}$ from $\cal I$. This process 
continues until $\cal I$ becomes empty. We will show that any vertex ordering 
$\prec$ which satisfies $(x,y) \in T \Rightarrow x \prec y$ for all vertices $x, y$ 
is a comparability ordering of $D$. 

\bigskip

\begin{breakablealgorithm}
\caption{}
\label{scdr}
\begin{algorithmic}[1]
    \Require
      A semicomplete digraph $D$.
    \Ensure
      Either a comparability ordering $\prec$ of $D$ or that
         ``$D$ is not a comparability digraph."
\If {$I=I^{-1}$ for any $I \in \mathcal{I}(D)$}
\State \Return ``$D$ is not a comparability digraph.''
\Else
\State initialize $T=\emptyset$ and
                  $\mathcal{I} = \{I \in \mathcal{I}(D):\ |I| > 1\}$.
\While {$\mathcal{I} \neq \emptyset$}
\If {$\exists\,(x,y), (y,z)\in T$ and $\exists I \in \mathcal{I}$ such that
           $(x,z)\in  I$}
\State choose $I$;
\Else \State arbitrarily choose $I \in \mathcal{I}$;
\EndIf
\State $T = T \cup I$ and
         $\mathcal{I}=\mathcal{I} \setminus \{I, I^{-1}\}$;
\EndWhile
\State \Return a vertex ordering $\prec$:\
           $(x,y) \in T \Rightarrow x \prec y$.
\EndIf
\end{algorithmic}
\end{breakablealgorithm}

\medskip   

\begin{thm} \label{correctness}
Algorithm \ref{scdr} correctly finds a comparability ordering of a semicomplete 
digraph $D$ for which $I \neq I^{-1}$ for all $I \in \mathcal{I}(D)$.
\end{thm}  

To prove Theorem \ref{correctness} we need a lemma. 
Let $C:\ (x_0,x_1), (x_1,x_2), \dots, (x_{k-2},x_{k-1}), (x_{k-1},x_0)$ 
be a circuit and $S$ be a set of pairs. A pair of $S$ is called {\em lonely} in 
$C$ if it is $(x_i,x_{i+1})$ for some $i$ and neither $(x_{i-1},x_i)$ nor 
$(x_{i+1},x_{i+2})$ is a pair of $S$. 

\begin{lem} \label{lem:structure}
Let $D$ be a semicomplete digraph and $I$ be a nontrivial implication class of $D$ 
with $I \neq I^{-1}$. Let $T$ be the union of some nontrivial implication classes of
$D$ disjoint from $I$ and $I^{-1}$. Suppose that $T$ contains no circuit but 
$T \cup I$ contains one. Then there exist three vertices such that 
$(a,b), (b,c)\in T$, $(a,c)\notin T$, and $[a,c]$ is nontrivial.
\end{lem}
\pf Let $C:\ (x_0,x_1), (x_1,x_2), \dots, (x_{k-2},x_{k-1}), (x_{k-1},x_0)$ be 
a shortest circuit in $T \cup I$. By assumption $T$ contains no circuit. Since 
$I \neq I^{-1}$, by Theorem \ref{thm:no-circuit-eq} $I$ contains no circuit.
Thus $C$ must contain a pair from $T \setminus I$ and a pair from 
$I\setminus T$. Since $T$ is disjoint from $I^{-1}$, $k \geq 3$. 

Suppose that for some $i$, $(x_i,x_{i+1})$ is a lonely pair of $T$ in $C$. Then 
$(x_{i-1},x_i), (x_{i+1},x_{i+2}) \in I$. Since $C$ is a shortest circuit in 
$T \cup I$, it has no chord in $T \cup I$, which implies that 
$(x_i,x_{i+2})\notin I$.
Applying the Triangle Lemma to $x_i, x_{i+1}, x_{i+2}$ with the fact
$(x_{i-1},x_i)\Gamma^* (x_{i+1},x_{i+2})$, we obtain that $(x_i,x_{i-1}) \in 
[x_i,x_{i+1}] \subseteq  T$. Thus $(x_i,x_{i-1})$ is in both $T$ and $I^{-1}$,
a contradiction to the fact that $T$ and $I^{-1}$ are disjoint.
Hence no pair of $T$ is lonely in $C$.

Suppose that for some $i$, $(x_i,x_{i+1}), (x_{i+1},x_{i+2})$ are pairs in $I$. 
By choosing different pairs if necessary we assume $(x_{i-1},x_i)$ is in $T$.
The choice of $C$ implies $(x_{i-1}, x_{i+1})\notin T\cup I$. 
Applying the Triangle Lemma to $x_{i-1}, x_i, x_{i+1}$ with the fact that
$(x_i,x_{i+1}) \Gamma^* (x_{i+1},x_{i+2})$, we obtain that 
$(x_{i-1},x_{i+1}) \in [x_{i-1},x_i] \subseteq T$, a contradiction.
Hence every pair of $I$ in $C$ is lonely.

Assume without loss of generality that $(x_{k-1},x_0) \in I$. Suppose that 
$(x_i,x_{i+1}) \in I$ is a pair in $C$ distinct from $(x_{k-1},x_0)$. Since every 
pair of $I$ in $C$ is lonely and no pair of $T$ in $C$ is lonely, 
$(x_0,x_1), (x_1,x_2), (x_{k-3},x_{k-2}), (x_{k-2},x_{k-1}) \in T$. 
So $2\leq i\leq k-4$. Applying the Triangle Lemma to $x_{k-2}, x_{k-1}, x_{0}$ with 
the fact that $(x_{k-1},x_{0})\Gamma^* (x_i,x_{i+1})$, we obtain that 
$(x_{k-2},x_i)\in [x_{k-2},x_{k-1}] \subseteq T$, a contradiction to the choice of 
$C$. Hence $(x_{k-1},x_0)$ is the only pair of $I$ in $C$. 

We claim that for each $0 \leq i \leq k-3$, $[x_{k-1},x_i] \neq [x_{k-1},x_{i+1}]$.
This is clearly true for $i = 0$ and $i = k-3$ for the choice of $C$. Suppose that 
for some $0 < i < k-3$, $[x_{k-1},x_{i-1}] \neq [x_{k-1},x_i] = [x_{k-1},x_{i+1}]$.
Applying the Triangle Lemma to $x_{i-1}, x_i, x_{k-1}$, we have 
$(x_{i-1},x_{i+1}) \in [x_{i-1},x_i] \subseteq T\cup I$, a contradiction to the
choice of $C$.    

Since $I$ is nontrivial, there exists a vertex $w$ such that 
$(x_{k-1},x_0) \Gamma (w,x_0)$ or $(x_{k-1},x_0) \Gamma (x_{k-1},w)$. We assume that
$(x_{k-1},x_0) \Gamma (w,x_0)$. (A similar proof applies when 
$(x_{k-1},x_0) \Gamma (x_{k-1},w)$.) Clearly, $w$ is not on $C$. 
Applying the Triangle Lemma to $x_1,x_0, x_{k-1}$ with the fact
$(x_0,w) \in [x_0,x_{k-1}]$, we get $(x_1,w) \in [x_1,x_{k-1}]$. 
Continuing this way with applications of the Triangle Lemma to 
$x_i,x_{i-1},x_{k-1}$, we get $(x_i,w) \in [x_i,x_{k-1}]$ for $i = 1, 2, \dots, k-3$.
In particular, $(x_{k-3},w) \in [x_{k-3},x_{k-1}]$. Therefore 
$x_{k-3},x_{k-2},x_{k-1}$ are three vertices such that 
$(x_{k-3},x_{k-2}), (x_{k-2},x_{k-1}) \in T$, $(x_{k-3},x_{k-1}) \notin T$, and 
$[x_{k-3},x_{k-1}]$ is nontrivial.
\qed

\medskip

{\bf Proof of Theorem \ref{correctness}}: 
We show first that $T$ obtained by Algorithm \ref{scdr} contains no circuit.
Suppose to the contrary that $T$ contains a circuit. Assume that $I$ is 
the first implication class added to $T$ that causes the occurrence 
of a circuit. That is, $T$ contains no circuit but $T\cup I$ contains 
one. By Lemma \ref{lem:structure}, there exist vertices $a, b, c$ such that
$(a,b), (b,c)\in T$, $(a,c)\notin T$, and $[a,c]$ is nontrivial. This means that
$I$ is chosen in steps 6-7 by Algorithm \ref{scdr}. Hence $I = [x,z]$ for some 
$(x,y),(y,z) \in T$.

Let $C:\ (x_0,x_1), (x_1,x_2), \dots, (x_{k-2},x_{k-1}), (x_{k-1},x_0)$ be a circuit
in $T\cup I$ that contains as few pairs in $I$ as possible. Assume with loss of 
generality that $(x_0,x_1) \in I$. Applying the Triangle Lemma to $y,x,z$ with
the fact that $(x_0,x_1) \in I = [x,z]$, we have $(x_0,y) \in [x,y] \subseteq T$
and $(y,x_1) \in [y,z] \subseteq T$. We see that 
$(x_0,y), (y,x_1), (x_1,x_2), \dots, (x_{k-2},x_{k-1}), (x_{k-1},x_0)$ is
a circuit in $T\cup I$ having fewer pairs in $I$ than $C$, a contradiction to the 
choice of $C$. Hence $T$ contains no circuit.

Since $T$ contains no circuit, a vertex ordering $\prec$ that satisfies
$(x,y) \in T \Rightarrow x \prec y$ exists. If $\prec$ is not a comparability
ordering of $D$, then for some $x \prec y \prec z$, $xy, yz \in A(D)$ and 
$xz \notin A(D)$ or $zy, yx \in A(D)$ and $zx \notin A(D)$. Either case implies that 
$(x,y) \Gamma (z,y)$. Thus $[x,y]$ and $[y,z]$ are nontrivial. 
Since $x \prec y \prec z$, $T$ contains both $(x,y)$ and $(y,z)$. Since 
$(x,y) \Gamma (z,y)$, $T$ contains $(z,y)$ as well. Hence $T$ contains both 
$(z,y)$ and $(y,z)$, contradictng the fact that $T$ contains no circuit.  
Therefore $\prec$ is a comparability ordering of $D$.
\qed

Combining Theorems \ref{thm:necessity}, \ref{thm:no-circuit-eq} and
\ref{correctness}, we obtain the following characterizations for semicomplete 
comparability digraphs, akin to Theorem~\ref{thm:undirected} for comparability 
graphs.

\begin{thm} \label{main}
The following statements are equivalent for a semicomplete digraph $D$:
\begin{enumerate}
\item $D$ is a comparability digraph;
\item the knotting graph $\tilde{K}_{D}$ is bipartite;
\item $D$ contains no odd closed knotting walk;
\item $I \neq I^{-1}$ for all $I\in \mathcal{I}(D)$;
\item $I$ contains no circuit for all $I\in \mathcal{I}(D)$.
\qed
\end{enumerate}
\end{thm}

We make some remarks on the complexity of Algorithm~\ref{scdr}. Assume that 
the input $D$ has $n$ vertices and is given by adjacency lists. For each 
$(x,y)\in Z_D$, we can find in $\mathcal{O}(n)$ time all pairs $(x,z)$ such that 
$(x,y)\Gamma (x,z)$. Thus we can compute all implication classes of $D$ in 
$\mathcal{O}(n^3)$ time. We separate the nontrivial implication classes from 
the trivial ones. By Theorem~\ref{thm:no-circuit-eq} an implication contains 
a circuit if and only if it contains a circuit of length 2. It follows that in 
$\mathcal{O}(n^2)$ time we can find out whether any implication class contains 
a circuit. We keep track of all pairs $(x,y) \in T$ for each vertex $x$ and for each 
update of $T$ and use it to determine whether $T$ contains $(x,y),(y,z)$ but not 
$(x,z)$ for some $x, y, z$. This can be done in $\mathcal{O}(n^3)$ time. Finally,
an $\mathcal{O}(n^2)$ time topological sort algorithm can be employed to find 
a vertex ordering $\prec$ of $D$ that satisfies $(x,y) \in T \Rightarrow x \prec y$.

\begin{thm} \label{complexity}
Algorithm \ref{scdr} can be implemented to run in $\mathcal{O}(n^3)$ time.
\qed
\end{thm}

\section{Concluding remarks and open problems}
\label{conclusion}

In this paper, we introduce the class of comparability digraphs as the digraph 
analogue of comparability graphs. We characterize comparability digraphs 
in terms of their knotting graphs. We extend the Triangle Lemma for graphs to
semicomplete digraphs. We use the Triangle Lemma for semicomplete digraphs to 
device an $\mathcal{O}(n^3)$ time recognition algorithm for semicomplete 
comparability digraphs. Finally we give a list of characterizations of semicomplete
comparability digraphs in the same way as for comparability graphs.   

Gallai's characterization of comparability graphs in terms of knotting graphs 
(see Theorem \ref{thm:undirected}) implies immediately a polynomial time recognition
algorithm for the class of graphs. Our characterization of comparability digraphs 
(see Theorem \ref{knotting-char}) is similar to Gallai's but it does not immediately
imply a polynomial time recognition algorithm for comparability digraphs.  
The digraph in Figure~\ref{noCs} suggests certain difficulties in finding one if 
it exists.

\begin{prob}
Determine the complexity of the recognition problem for comparability digraphs.
\end{prob}

Gallai \cite{gallai} has also given a forbidden subgraph characterization of 
comparability graphs. The minimal forbidden subgraphs for comparability graphs 
comprise 19 types of graphs (cf. \cite{gallai}). We propose as open problems 
the following:

\begin{prob}
Find all minimal forbidden subdigraphs of (semicomplete) comparability digraphs.
\end{prob}

\end{CJK}

\begin{thebibliography}{99}

\bibitem{afarber} R.P. Anstee and M. Farber, Characterizations of totally
balanced matrices,
{\it J. Algorithms} 5 (1984) 215 - 230.

\bibitem{bdr} D. Bechet, P. de Groote, and C. Retor\'e, A complete axiomatisation 
of the inclusion of series-parallel partial orders, 
{\it Rewriting Techniques and Applications}, {\it LNCS} 1232 (1997) 230 - 240. 

\bibitem{boe} D. Boechner, Oriented threshold graphs,
{\it Australasian J. Combinatorics} 71 (2018) 43 - 53.

\bibitem{dein} V. Deineko, R. Rudolf, and G.J. Woeginger,
A general approach to avoiding two by two submatrices,
{\it Computing} 52 (1994) 371 - 388.

\bibitem{gallai} T. Gallai, Transitiv orientierbare Graphen, 
{\it Acta Mathematica Academiae Scientiarum Hungarica} 18 (1967) 25 - 66.

\bibitem{golumbic77} M.C. Golumbic, The complexity of comparability graph 
recognition and coloring, {\it Computing} 18 (1977) 199 - 208.

\bibitem{golumbic} M.C. Golumbic, {\bf Algorithmic Graph Theory and
Perfect Graphs}, Academic Press (1980).

\bibitem{chordal-di} L. Haskins, D.J. Rose, Toward characterization of perfect 
elimination digraphs, {\it SIAM J. Computing} 2 (1973) 217 - 224.

\bibitem{recg-lexico} P. Hell and J. Huang, Lexicographic orientation and 
representation algorithms for comparability graphs, proper circular arc graphs, and 
proper interval graphs, {\it J. Graph Theory} 20 (1995) 361 - 374.

\bibitem{hhlim} P. Hell, J. Huang, J.C.-H. Lin, and R.M. McConnell, 
Bipartite analogues of comparability and cocomparability graphs, 
{\it SIAM J. Discrete Math.} 34 (2020) 1969 - 1983.

\bibitem{minorder} P. Hell, J. Huang, R.M. McConnell, and A. Rafiey, 
Min-orderable digraphs, 
{\it SIAM J. Discrete Math.} 34 (2020) 1710 - 1724.

\bibitem{hks} A.J. Hoffman, A.W. Kolen, and M. Sakarovitch, Totally-balanced
and greedy matrices,
{\it SIAM J. Algebraic and Discrete Methods} 6 (1985) 721 - 730.

\bibitem{krw} B. Klinz, R. Rudolf, and G.J. Woeginger,
Permuting matrices to avoid forbidden submatrices,
{\it Discrete Applied Math.} 60 (1995) 223 - 248.

\bibitem{lamar} M.D. Lamar, Split digraphs,
{\em Discrete Math.} 312 (2012) 1314 - 1325.

\bibitem{lss} M.C. Lin, F.J. Soulignac, and J.L. Szwarcfiter,
Arboricity, h-index, and dynamic algorithms,
{\it Theor. Compt. Sci.} 426/427 (2012) 75 - 90.

\bibitem{lubiw1} A. Lubiw, $\Gamma$-free matrices,
Master's thesis, University of Waterloo, 1982.

\bibitem{1999Modular} R.M. McConnell and J.P. Spinrad, Modular decomposition and
transitive orientation, {\it Discrete Math.} 201 (1999) 189 - 241.

\bibitem{sdrw} M. Sen, S. Das, A.B. Roy, and D.B. West, Interval digraphs:
An analogue of interval graphs, {\it J. Graph Theory} 13 (1989) 189 - 202.

\bibitem{Circular-arc-di} M. Sen, S. Das, and D.B. West, Circular-arc digraphs: 
a characterization, {\it J. Graph Theory}, 13 (1989) 581 - 592.

\end{thebibliography}
\end{document}